\documentclass{article}

\usepackage[T1]{fontenc}
\usepackage[english]{babel}
\usepackage{amssymb}
\usepackage{amsmath}
\usepackage{latexsym}
\usepackage{marvosym}
\usepackage{amsfonts}
\usepackage{graphicx}
\usepackage{color}
\usepackage{fullpage}
\usepackage[latin1]{inputenc}
\usepackage{eurosym}
\usepackage{enumerate}
\usepackage{fancyhdr}
\usepackage{bbold}
\usepackage{tikz}
\usepackage{titletoc}
\usepackage{natbib}

\usepackage[pdftex=true,hyperindex=true,colorlinks=true,linkcolor=blue,bookmarks
=true]{hyperref}

\definecolor{fgblue}{rgb}{0.2, 0.3, 0.6}

\hypersetup{
backref=true,
pagebackref=true,
hyperindex=true,
colorlinks=true,
breaklinks=true, 
urlcolor= fgblue,
linkcolor= fgblue,
citecolor=fgblue,
bookmarks=true,
bookmarksopen=true,
}

\usepackage{lmodern}

\usepackage{eulervm}

\newtheorem{definition}{D\'efinition}[section]
\newtheorem{proposition}[definition]{Proposition}
\newtheorem{theorem}[definition]{Theorem}
\newtheorem{corollary}[definition]{Corollary}
\newtheorem{hypothese}[definition]{Assumption}

\newtheorem{remarque}[definition]{Remark}

\newtheorem{lemma}[definition]{Lemma}

\newcommand{\N}{\mathbb{N}}
\newcommand{\E}{\mathbb{E}}
\newcommand{\esp}{\mathbb{E}}
\newcommand{\R}{\mathbb{R}}
\newcommand{\prob}{\mathbb{P}}

\newcommand{\e}{\varepsilon}

\newcommand{\proof}{\noindent\textit{\textbf{Proof.}}~}

\title{Semi-parametric inference for the absorption\\
features of a growth-fragmentation model}
\author{Romain Aza\"\i{}s$^\dagger$\footnote{Corresponding author. E-mail address: \url{romain.azais@gmail.com}}~~and Alexandre Genadot$^\ddagger$}
\date{\small $\dagger$\,Inria Sophia Antipolis M\'editerran\'ee, Team Virtual Plants\\
95, Rue de la Gal\'era, 34 095 Montpellier\\
\vspace{2mm}
$\ddagger$\,Laboratoire de Probabilit\'es et Mod\`eles Al\'eatoires, Universit\'e Pierre et Marie Curie, Paris 6\\
4, Place Jussieu, 75 005 Paris}

\begin{document}

\maketitle

\vspace{-0.5cm}

\begin{abstract}
In the present paper, we focus on semi-parametric methods for estimating the absorption probability and the distribution of the absorbing time of a growth-fragmentation model observed within a long time interval. We establish that the absorption probability is the unique solution in an appropriate space of a Fredholm equation of the second kind whose parameters are unknown. We estimate this important characteristic of the underlying process by solving numerically the estimated Fredholm equation. Even if the study has been conducted for a particular model, our method is quite general.
\end{abstract}

%\vspace{1mm}



%\vspace{-0.4cm}

%\dottedcontents{section}[3.8em]{\addvspace{0.8pt}}{2.3em}{1em}
%\addtocontents{toc}{\protect\setcounter{tocdepth}{1}}
%\tableofcontents

\section{Introduction}

The present paper is dedicated to a statistical approach for a particular growth-fragmentation model $(X_t)_{t\geq0}$ for which $\Gamma=[0,1]$ is an absorbing set. The motion of $(X_t)_{t\geq0}$ involves expontential growth and {downward} random jumps at random times. From the observation of only one trajectory of the process starting from $X_0=x\notin\Gamma$, we propose to estimate the probability of absorption $p(x)$ and the distribution $(t_m(x))_{m\geq1}$ of the hitting time of $\Gamma$. {Since the process $(X_t)_{t\geq0}$ grows exponentially between two consecutive downward jumps, absorption may only occur at jumps,} and $t_m(x)$ is the probability of absorption at the $m^{\text{th}}$ jump of the process $(X_t)_{t\geq0}$ starting from $X_0=x$. This stochastic process has already been introduced by \cite{TJRI} as a theoretical model for insurance, for which $\Gamma$ is called \textit{area of poverty}. \cite{TJRI} establish that the ruin probability $p(x)$ is the solution of an integral equation that they solve numerically. In many aspects, our work and the aforementioned paper are different and complementary. Indeed, we adopt a statistical approach: from the observation of the process within a long time, we propose a semi-parametric procedure for estimating the two main characteristics $p(x)$ and $(t_m(x))_{m\geq1}$ of this model. In addition, we would like to highlight that naive estimates obtained from the empirical versions of these features would not have worked in the framework of a long time observation.

The growth-fragmentation model that we consider is a particular case of piecewise-deterministic Markov processes (PDMP's). The PDMP's were first introduced in the literature by \cite{DAVIS} in the eighties as a general class of non-diffusion stochastic models. They form a family of continuous-time Markov processes involving deterministic motion punctuated by random jumps, which occur either in a Poisson-like fashion with nonhomogeneous rate or when the deterministic flow hits the boundary of the state space. The motion of a PDMP depends on three local characteristics namely the jump rate, the flow and the transition kernel. A suitable choice of the state space and these three features provides a large number of stochastic models covering various applications, for example in reliability \citep{ChiquetLimnios,DeS} or in neurosciences \citep{buckwar2011exact,genadot2012averaging,riedler2012limit}. These processes have been heavily studied both from theoretical and applied perspectives. One may refer the reader to the recent papers \citep{benaim2014,cloezhairer,costa:siam} about their ergodicity properties and to \citep{brandejsky:spa,CD:book} for references on optimal control and stopping for this class of stochastic models. There are only a few papers which deal with statistical methods for PDMP's. Without attempting to make an exhaustive survey about these methods, we refer the reader to \citep{AzaisESAIM14,AzaisSJS14,Doumic11,Doumic12,JACOB} and the references therein. \cite{Doumic11,Doumic12} focus on estimation procedures for ergodic size-structured models, while \cite{AzaisESAIM14} and \cite{AzaisSJS14} are interested in the nonparametric estimation of the transition measure and of the conditional distribution of the inter-jumping times for a general PDMP observed within a long time, under some ergodicity conditions. We would like to emphasize that these methods should not work in our non-ergodic framework. \cite{JACOB} computes likelihood processes for observation of PDMP's without boundary jumps. This approach could lead to some estimation methods in the parametric case.

Our framework is semi-parametric. Indeed, the features of the PDMP that we consider are defined from a density function $G$ on $[0,1]$ and a parameter $\lambda>0$. Our approach consists in building estimators of the probability of absorption and of the distribution of the absorbing time from some estimators of the features $G$ and $\lambda$. %We will give an example of such estimators of $G$ and $\lambda$ which satisfy the few asymptotic conditions that we impose.
For this purpose, we show that the probability of absorption is solution of a Fredholm equation of the second kind that we propose to estimate. This allows us to build our estimator of the probability of absorption $p(x)$ as the numerical solution of this estimated equation. Thus, our estimator $\widehat{p}_{n,m}(x)$ depends on both the number $n$ of available data and the number $m$ of iterations of the numerical resolution algorithm. This estimator is defined in $(\ref{eq:pnm})$. Therefore, the error is the combination of a statistical error and a numerical error. The result of convergence of $\widehat{p}_{n,m}(x)$ towards $p(x)$ is stated in Theorem \ref{prop:estim:r}. Of course, the absorption probability is an important feature of the model considered, but the law of the hitting time of $\Gamma$ is a significant complementary information. We also build an estimator of $t_m(x)$, the probability for the process starting from $X_0=x$, to be absorbed at the m$^{th}$ jump. This estimator, $\widehat{t}_{m,n}(x)$, is defined in $(\ref{eq:tnm})$ and the result of convergence of $\widehat{t}_{m,n}(x)$ towards $t_m(x)$ is stated in Theorem \ref{prop:estim:t}. To estimate the two main quantities of interest $p(x)$ and $(t_m(x))_{m\geq1}$, we need to estimate at first the post-jump location transition kernel of the process $(X_t)_{t\geq0}$. This estimator is defined in $(\ref{eq:estim:R})$ and convergence results are stated in Proposition \ref{prop:maj:R}.

We take care to illustrate our theoretical results by numerical simulations. For a given example, we carry out the whole estimation procedure for the quantities of interest $p(x)$ and $(t_m(x))_{m\geq1}$, on a hundred replications and several numbers of observed data ($50$, $75$ and $100$ observed jumps). We present the corresponding numerical results in Section \ref{sec:simulations}. {Let us notice that the numerical estimates are close to the expected values despite the low number of data}, especially with respect to the sample sizes used in the works \citep{AzaisESAIM14,AzaisSJS14,Doumic11,Doumic12}.

The paper is organized as follows. We begin in Section \ref{s:problem} with the problem formulation. We present the absorbing growth-fragmentation model that we focus on in Subsection \ref{ss:model} and the semi-parametric framework that we choose in Subsection \ref{ss:semiparam}. In addition, two applications of the present stochastic model are presented in Subsection \ref{subs:appl}. Section \ref{sec:main:res} is devoted to the presentation of our main results. The results about the estimation of the transition kernel of the post-jump locations are presented in Subsection \ref{ss:R}. Subsections \ref{subs:abs} and \ref{subs:hit} gather the two main results of the paper: the construction and convergence of our estimators for the probability and the time of absorption. Our results are illustrated by numerical simulations within a particular setting in Section \ref{sec:simulations}. Some additional results and the proofs of the main results have been deferred {until} Appendix \ref{sec:annexe}, \ref{app:proofs} and \ref{sec:discu}.

\section{Problem formulation}
\label{s:problem}

In this section, we present the absorbing growth-fragmentation model that we focus on and some applications. In addition, we present the semi-parametric framework that we choose for investigating the statistical inference for this process observed within a long time interval.

\subsection{A growth-fragmentation model}
\label{ss:model}

The growth-fragmentation model that we consider is the Markov process $(X_t)_{t\geq0}$ defined on some pro\-ba\-bility space $(\Omega, \mathcal{F},\prob)$ -- from a pro\-ba\-bility density function $G$ on the interval $[0,1]$ and two real numbers $r, \lambda>0$ -- by its extended generator $\mathcal{A}$ as follows,
%%%
\begin{equation}\label{eq:generator}
\mathcal{A}f(x)=r(x-1)^+f'(x)+\lambda\int_0^1\big(f(zx)-f(z)\big)\,G(z)\,dz
\end{equation}
%%%
with $\xi^+=\xi\vee0$, and for any smooth function $f\,:\,\R\to\R$ (see \citep{DAVIS} for full details on the domain of $\mathcal{A}$). {The generator} $(\ref{eq:generator})$ {describes a process having an exponential growth at rate $r$ and which can only decrease by means of downward jumps. The downward jumps occur at rate $\lambda$ and instantaneously reduce the process of a certain percentage distributed according to the density $G$. The positive part present in equation} $(\ref{eq:generator})$ {implies that the exponential growth is annihilated below the threshold $1$, preventing the process to escape the domain $[0,1]$.} Some applications of this stochastic model are presented in Subsection \ref{subs:appl}. Let us notice that this process has already been introduced as a stochastic model for insurance in \citep{TJRI}. In the sequel, we propose to describe the motion of our growth-fragmentation model as the dynamic of a piecewise-deterministic Markov process (PDMP), see the book \citep{DAVIS} and the references therein.

In most cases, the dynamic of a one-dimensional real-valued PDMP is described by its three local features $(\widetilde{\lambda},\mathcal{Q},\Phi)$.
\begin{itemize}
	\item $\Phi : \R\times\R_+\to\R$ is the deterministic flow. It satisfies,
		$$\forall\,\xi\in\R,~\forall\,s,t \geq 0,~ {\Phi(\xi,0)=\xi \quad\text{and}\quad}  \Phi(\xi,t+s) = \Phi(   \Phi(\xi,t) , s ).$$
		%In particular, $\Phi(\xi,0)=\xi$.
		%For each $\xi\in E$, $t^\star(\xi)$ denotes the deterministic exit time from $E$:
		%$$ t^\star(\xi) = \inf \{t>0~:~\Phi(\xi,t) \in \partial E\},$$
		%with the usual convention $\inf\emptyset = +\infty$.
	\item $\widetilde{\lambda}: \R\to\R_+$ is the jump rate. {It satisfies},
		$$\forall\,\xi\in\R,~\exists\,\varepsilon>0,~ \int_0^\varepsilon \widetilde{\lambda}\big( \Phi(\xi,s) \big) ds < \infty .$$
	\item $\mathcal{Q}$ is a Markov kernel on $\R$ which satisfies,
		$$\forall\,\xi\in\R,~\mathcal{Q}(\xi,\{\xi\})=0.$$
\end{itemize}
Starting from $X_0=x$, the motion can be described as follows.
The first jump time $T_1$ is a positive random variable whose survival function is,
$$ \forall\,t \geq 0, ~ \prob( T_1 > t \,|\, X_0=x) = \exp \left(  - \int_0^t \widetilde{\lambda}( \Phi(x,s) )d s \right).$$
This jump time occurs in a Poisson-like fashion with nonhomogeneous rate $\widetilde{\lambda}$.
One chooses a real-valued random variable $Z_1$ according to the distribution $\mathcal{Q}(\Phi(x,T_1),\cdot)$.
Let us remark that the post-jump location $Z_1$ depends on the interarrival time $T_1$,
via the deterministic flow starting from $X_0=x$. The trajectory between the times $0$ and $T_1$ is given by
\begin{displaymath}
X_t = \left\{ 
\begin{array}{lll}
\Phi(x,t) 	& \text{for}& 0\leq t < T_1, \\
Z_1		& \text{for}& t=T_1.
\end{array}
\right.
\end{displaymath}
Now, starting from $X_{T_1}$, one may choose the interarrival time $S_2 = T_2-T_1$ and the post-jump location $Z_2$ in a similar way as before,
and so on. The randomness of such a process is only given by the jump mechanism.

In our particular case, one may easily compute from $(\ref{eq:generator})$ the local features of the growth-fragmentation model $(X_t)_{t\geq0}$ that we consider. They are given by
\begin{equation}\label{eq:def:flow}
{\Phi(x,t)} =
\left\{
\begin{array}{lll}
(x-1) \exp(rt)+1	& \text{if}			&x>1,\\
x			&	\text{else,}	&
\end{array}
\right.
\end{equation}
\begin{equation}\label{eq:def:Q}
\widetilde{\lambda}(x)=\lambda\qquad\text{and}\qquad\mathcal{Q}(x,dy)=\frac{1}{x}G\left(\frac{y}{x}\right)dy.
\end{equation}

\begin{figure}[!h]
\centering
\includegraphics[height=5cm,width=11cm]{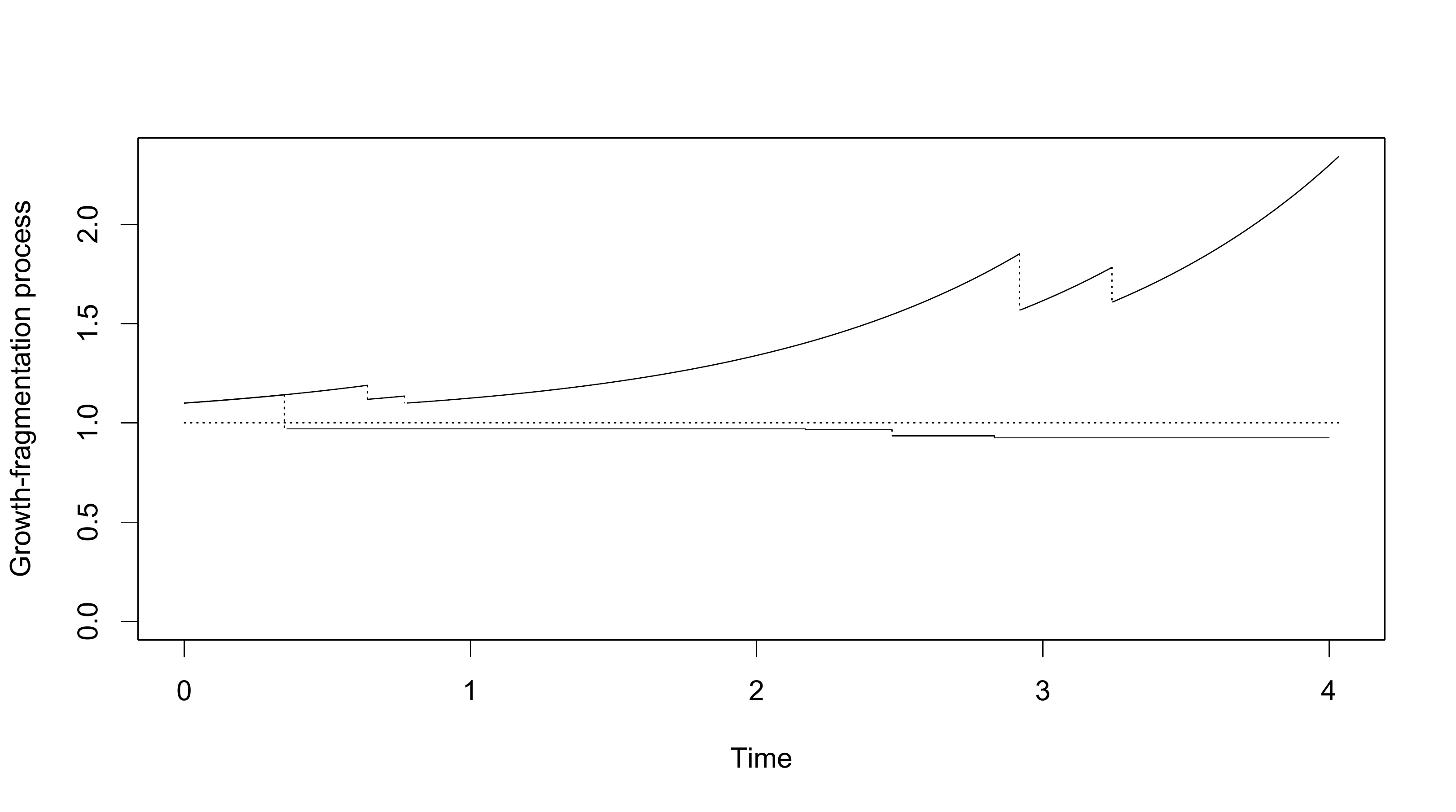}
\caption{Two trajectories for the growth-fragmentation model {starting from $X_0=1.1$ and whose characteristics are $\lambda=r=1$ and $G(u)=11u^{10}$, as in the numerical example investigated in Section} \ref{sec:simulations}. One of them is absorbed at the first jump time, while the other one seems to escape the trapping set.}
\label{fig:ex}
\end{figure}

\noindent
{Notice that the flow $\Phi$ given by equation} $(\ref{eq:def:flow})$ {arises, according to} $(\ref{eq:generator})${, as the solution of the following ordinary differential equation,
$$
y'_t=r(y_t-1)^+\quad\text{and}\quad y_0=x,
$$
describing the kinetic of the process between jumps. Moreover, we see} in particular from (\ref{eq:def:Q}), {that} the rate {of jump} $\widetilde{\lambda}$ is homogeneous {and that jumps are downward since $G$ is a density on $[0,1]$}. Therefore, the sequence of the interarrival times $(S_n)_{n\geq1}$ is independent and exponentially distributed with rate $\lambda$. In addition, the particular form of the transition kernel $\mathcal{Q}$ $(\ref{eq:def:Q})$ implies that the sequence of the random loss fractions $(Y_n)_{n\geq1}$ defined from,
\begin{equation}\label{eq:randomlossfraction}
\forall\,n\geq1,~Z_n={\Phi(Z_{n-1},S_n)}\,Y_n,
\end{equation}
is independent and independent of $(S_n)_{n\geq1}$ with common distribution $G$ {(see Lemma} \ref{lem:Yn:iid}{)}. As a consequence, the dynamic of the PDMP $(X_t)_{t\geq0}$ may be summarized by the observation of the independent sequences $(S_n)_{n\geq1}$ and $(Y_n)_{n\geq1}$. Two possible trajectories of the process $(X_t)_{t\geq0}$ are given in Figure \ref{fig:ex}.

%\begin{figure}[!h]
%\centering
%\includegraphics[height=5cm,width=11cm]{01_trajectoires.pdf}
%\caption{Two trajectories for the growth-fragmentation model. One of them is absorbed at the first jump time, while the other one seems to escape the trapping set.}
%\label{fig:ex}
%\end{figure}

This continuous-time Markov process is called absorbing because the motion may reach the absorbing interval $\Gamma=[0,1]$ from any initial value $X_0=x$. In this paper, we focus on the estimation of the absorption probability (see Subsection \ref{subs:abs}) and of the distribution of the hitting time of $\Gamma$ (see Subsection \ref{subs:hit}) from the observation of only one trajectory of the process within a long time, that is to say from the observation of the sequences $(S_n)_{n\geq1}$ and $(Y_n)_{n\geq1}$. {We would like to emphasize that this is far from obvious to estimate these absorption features from only one trajectory of the process, which has been either absorbed or not.}

\begin{remarque}
{The growth-fragmentation process described in this section satisfies the following Poisson driven stochastic differential equation,
\begin{equation}\label{eq_Poisson}
dX_t=r(X_t-1)^+dt+X_{t^-}\int_0^1(z-1)\,\mathcal{N}(dt\times dz), \quad t\geq0,
\end{equation}
with initial condition $X_0=x$. The measure $\mathcal{N}$ is a Poisson random measure on $\R_+\times [0,1]$ with intensity $\nu$, the measure on $\R_+\times [0,1]$ defined by
$$
\nu(dt\times dz)=\lambda dt\,G(z)dz.
$$
Equation $(\ref{eq_Poisson})$ provides another dynamical representation for the process $(X_t)_{t\geq0}$ which may be used, for example, for simulation purposes. Even a closed, but not very tractable, form for $(X_t)_{t\geq0}$ may be derived from $(\ref{eq_Poisson})$. This will not be used in the present paper.}
\end{remarque}

\subsection{Semi-parametric framework}
\label{ss:semiparam}

In all the sequel, we assume that we observe a PDMP $(X_t)_{t\geq0}$ within a long time interval. In other words, we observe the $n$ first terms of the sequence of the interarrival times $(S_k)_{k\geq1}$ and of the sequence of the independent random loss fractions $(Y_k)_{k\geq1}$ defined by $(\ref{eq:randomlossfraction})$.

From these independent observations, we propose to estimate the features $G$ and $\lambda$. In the rest of the paper, we consider an estimator $\widehat{G}_n$ of the density function $G$ and an estimator $\widehat{\lambda}_n$ of the rate $\lambda$, computed from the $n$ first loss events, that is to say from $S_1,\,\dots,\,S_n$ and $Y_1,\,\dots,\,Y_n$. In some of our convergence results, we impose a few conditions on the asymptotic behaviors of $\widehat{G}_n$ and $\widehat{\lambda}_n$.
When one of the assumptions {is used}, this {is specified} in the statement of the result.

\vspace{1mm}

\begin{minipage}[l]{0.38\textwidth}
\begin{description}
\item [$(C_1^\lambda)$] $\widehat{\lambda}_n\in[\lambda_\ast,\lambda^\ast]$, with $\lambda_\ast>0$.
\item [$(C_2^\lambda)$] ${\widehat{\lambda}_n-\lambda}$ tends to $0$ in probability.
%\item [$(C_3^\lambda)$] $\widehat{\lambda}_n$ almost surely tends to $0$.
\end{description}
\end{minipage}
\begin{minipage}[l]{0.62\textwidth}
\begin{description}
\item [$(C_1^G)$] $\|\widehat{G}_n-G\|_\infty$ tends to $0$ in probability.
\item [$(C_2^G)$] $\int_0^1 |G(u)-\widehat{G}_n(u)|\,u^{-1}\,du$ tends to $0$ in probability.
%\item [$(C_2^G)$] $\|\widehat{G}_n-G\|_\infty$ almost surely tends to $0$.
\end{description}
\end{minipage}

\begin{remarque}
{The condition $(C^\lambda_1)$ states that the sequence of the estimators $\widehat{\lambda}_n$ is uniformly bounded in $(0,\infty)$. The three other conditions, $(C^\lambda_2)$ and $(C^G_{1,2})$, concern the convergence in probability of the estimators $\widehat{\lambda}_n$ and $\widehat{G}_n$. Regarding $\widehat{G}_n$, this convergence in probability is stated with respect to two different norms: the sup-norm for $(C^G_1)$ and a weighted $L^1$-norm for $(C^G_2)$. This latter convergence mode implies that the probability to jump back close to zero should not be too large.}
\end{remarque}

\noindent
In the present paper, we are not interested in the demonstration of the asymptotic properties of the estimates $\widehat{G}_n$ and $\widehat{\lambda}_n$ but in the estimation of some characteristics of the PDMP $(X_t)_{t\geq0}$ from these estimates. In particular, we establish in Theorems \ref{prop:estim:r} and \ref{prop:estim:t} that the convergences in probability of $\widehat{G}_n$ and $\widehat{\lambda}_n$ may be transferred to our estimators of the absorption probability and of the distribution of the hitting time of $\Gamma$. Thus, we do not investigate the properties or the choice of the estimators of $G$ and $\lambda$. Nevertheless, the assumptions that we impose on $\widehat{G}_n$ and $\widehat{\lambda}_n$ are non restrictive {as we proceed to show, see also Appendix \ref{sec:discu}}. For the sake of readability, we introduce the following notations: $\overline{S}_n$ denotes the empirical mean of the $n$ first interarrival times, while the projection $\pi_{[a,b]}(x)$ is defined by
$$\pi_{[a,b]}(x) =
\left\{
\begin{array}{clc}
a & \text{if} & x<a,\\
b & \text{if} & x>b,\\
x & \text{else}. & 
\end{array}
\right.
$$
In the case where we know two real numbers $\lambda_\ast$ and $\lambda^\ast$ such that $0<\lambda_\ast<\lambda<\lambda^\ast$, the truncated maximum likelihood estimator given by
\begin{equation}\label{eq:tmle}\widehat{\lambda}_n^{ml} = \pi_{[\lambda_\ast,\lambda^\ast]}\big(\hspace{0.5mm}\overline{S}_n^{-1}\hspace{0.2mm}\big),\end{equation}
obviously satisfies both the conditions $(C_{1,2}^\lambda)$. Furthermore, the Parzen-Rosenblatt estimator of a uniformly continuous density satisfies the assumption $(C_{1}^G)$ whenever the bandwidth $(h_n)_{n\geq1}$ is such that
$$\sum_{n\geq1} \exp\left(-\delta n h_n^2\right) < \infty,$$
for any $\delta>0$ (see \citep{Wied2010Consistency} for instance). Finally, we will establish in Appendix \ref{sec:discu} that the convergence $(C_2^G)$ is also satisfied for the Parzen-Rosenblatt estimator under some additional conditions (see Assumption \ref{hyp:G2}) on the density of interest $G$.

%%%%%%%%%%%%%%%%%%%%%%%%%
\begin{remarque}
{The situation that we consider in the present paper may be generalized to estimate the absorption features of a large number of stochastic processes. Indeed, our method may be applied to any absorbing jump process satisfying both the following conditions:
\begin{itemize}
\item the absorption occurs only at (downward or upward) jumps;
\item the motion of the process depends on an easily estimable quantity $\Theta$ in a parametric, semi-parametric or nonparametric setting. In our semi-parametric framework, $\Theta=(\lambda,G)$ is estimated from the independent random variables $(S_i)_{1\leq i\leq n}$ and $(Y_i)_{1\leq i\leq n}$. Let us denote $\widehat{\Theta}_n$ an estimator of $\Theta$.
\end{itemize}
Indeed, we establish in Proposition \ref{prop:eq:int} that the absorption probability is solution of an integral equation without using the specific kinetic of the model but only the fact that the absorption occurs at jumps. This solution depends on the transition kernel of the embedded chain of the jump process, and thus on the quantity $\Theta$ which governs the dynamic. As a consequence, one may estimate the absorption probability by plugging the estimate $\widehat{\Theta}_n$ in the expression of this solution. Nevertheless, the transfer of the asymptotic properties of $\widehat{\Theta}_n$ to the estimator of the absorption features strongly depends on the specific motion of the chosen model.}
\end{remarque}
%%%%%%%%%%%%%%%%%%%%%%%%%

\subsection{Some applications}
\label{subs:appl}

In this part, we present two applications of the growth-fragmentation model that we focus on in the present paper.

\subsubsection{A ruin theoretical model}\label{sec_ex_eco}

The growth-fragmentation process that we consider has been introduced by \cite{TJRI} for modeling a capital subject to random heavy loss events. We do not attempt to give an exhaustive survey about this model, but refer the reader to the paper \citep{TJRI} and the references therein. We consider an individual household whose income $I_t$ at time $t$ may be split into
\begin{equation}\label{eq:construction:1}
I_t = C_t + {\sigma_t}
\end{equation}
where $C_t$ denotes the consumption and ${\sigma_t}$ is the savings.
Consumption is assumed to evolve according to
\begin{equation}\label{eq:construction:2}
C_t=
\left\{
\begin{array}{lll}
I_t & \text{if}&I_t\leq I^\ast\\
I^\ast + a(I_t-I^\ast) & \text{else,}&
\end{array}
\right.
\end{equation}
where $I^\ast$ is the critical income level and $0<a<1$.
If the income is smaller than $I^\ast$, the whole income is used for consumption.
We denote by $X_t$ the accumulated capital up to time $t$.
The capital evolves according to
\begin{equation}\label{eq:construction:3}
\frac{dX_t}{dt} = c {\sigma_t},
\end{equation}
where $0<c<1$, while the income evolves with
\begin{equation}\label{eq:construction:4}
I_t = b X_t,
\end{equation}
with $b>0$. Finally, from $(\ref{eq:construction:1})$, $(\ref{eq:construction:2})$, $(\ref{eq:construction:3})$ and $(\ref{eq:construction:4})$, the capital $X_t$ satisfies the ordinary differential equation,
\begin{equation}\label{eq:construction:ode}
\frac{dX_t}{dt} = r (X_t - x^\ast)^+,
\end{equation}
where $x^\ast=I^\ast/b$ and $r=(1-a) b c$.
Now, we assume that the capital $X_t$ is subject to catastrophic events,
which occur in a Poisson-like fashion with homogeneous rate $\lambda$.
When an event occurs at time $t$, the capital $X_t$ is reduced by a random fraction
whose distribution is described by its probability density function $G$.
After the loss event, the process starts again according to $(\ref{eq:construction:ode})$. We obtain a PDMP for which the interval $[0,x^\ast]$ is an absorbing set, called \textit{area of poverty}. Indeed, once the process is below the critical capital $x^\ast$, the next events will reduce the capital and the process will never again reach values above $x^\ast$.

\subsubsection{A Malthusian evolution}\label{sec_ex_bio}

Here is an example from population dynamics. Consider a population whose total number of individuals at time $t$ is $X_t$. We assume that there exists a certain extinction threshold $x^*$ below which the population will almost surely extinct. One of the simplest dynamic for population models is that of Malthus (see \citep{MURRAY} for instance). In this model, a super-threshold population grows exponentially according to the equation $(\ref{eq:construction:ode})$.
%$$
%\frac{dX_t}{dt}=r(X_t-x^*)^+.
%$$
Of course, if you not {enrich} the model with new assumptions, starting from a super-threshold population, the population shall never be extinguished. However, note that this simple exponential growth model describes pretty well the growth of the human population over the past centuries. One way to make the model more realistic is to assume that disasters can happen. If we think {of} the human population, such disasters may correspond to epidemics, famines or wars. As above for the ruin theoretical model, we consider that these disasters occur in a Poisson-like fashion and that when a disaster happens, it has the effect of instantaneously {reducing} the population of a certain proportion. In this case, it may happen that starting from a super-threshold population, at some point in the time, a disaster affects the population so much so that it falls below the extinction threshold and therefore that extinction occurs.

\section{Main results}\label{sec:main:res}

We present here our main results on the estimation of the absorption probability and of the hitting time of $\Gamma$ for the PDMP $(X_t)_{t\geq0}$. First, we focus on a procedure for estimating the Markov kernel of the post-jump locations of this process. Almost all the proofs have been deferred until Appendix \ref{app:proofs}.

\subsection{Transition density of the post-jump locations}
\label{ss:R}

We are interested in the estimation of the transition kernel $\mathcal{R}$ of the Markov chain of the post-jump locations $(Z_n)_{n\geq0}$. Recall that the sequence of the jump times of the PDMP $(X_t)_{t\geq0}$ is $(T_n)_{n\geq0}$ and the post-jump locations are defined, for any $n\geq0$, by $Z_n=X_{T_n}$. The Markov kernel $\mathcal{R}$ of $(Z_n)_{n\geq0}$ is given, for any $x\in\R_+^\ast$ and $A\in\mathcal{B}(\R_+^\ast)$, by
$$\mathcal{R}(x,A) = \prob(Z_{n+1}\in A\,|\,Z_n=x).$$
In the rest of the paper, we impose the following condition on the probability density $G$. {This assumption means that the probability to jump back close to zero should not be too large.}
\begin{hypothese}\label{hyp:G}
We assume that the density $G$ is bounded and such that $\int_0^1 G(u)\,u^{-1} du<\infty$.
\end{hypothese}

\noindent
First, we show that $\mathcal{R}$ may be directly computed from the characteristics $G$ and $\lambda$.

\begin{proposition}\label{prop:def:ker}
The transition kernel of the Markov chain $(Z_n)_{n\geq0}$ satisfies $\mathcal{R}(x,dy)=\mathcal{R}(x,y)dy$, with
\begin{equation}\label{eq:calcul:R}
\mathcal{R}(x,y) =
\left\{
\begin{array}{lll}
\frac{1}{x}G\left(\frac{y}{x}\right) & \text{if}&x\leq1,\\
\frac{\lambda}{r}(x-1)^{\lambda/r}\,\left[ \int_{0}^{y/x\wedge 1} G(u) u^{\lambda/r}\,(y-u)^{-\lambda/r-1}du\right]
& \text{else.}&
\end{array}
\right.
\end{equation}
\end{proposition}

\noindent
Therefore, by virtue of $(\ref{eq:calcul:R})$, one may propose to estimate the transition density $\mathcal{R}(x,y)$ by
\begin{equation}\label{eq:estim:R}
\widehat{\mathcal{R}}_n(x,y) =
\left\{
\begin{array}{lll}
\frac{1}{x}\widehat{G}_n\left(\frac{y}{x}\right) & \text{if}&x\leq1,\\
 \frac{\widehat{\lambda}_n}{r}(x-1)^{\widehat{\lambda}_n/r}\,\int_{0}^{y/x\wedge 1} \widehat{G}_n(u) u^{\widehat{\lambda}_n/r}\,(y-u)^{-\widehat{\lambda}_n/r-1}du
& \text{else,}&
\end{array}
\right.
\end{equation}
where $\widehat{G}_n(\xi)$ and $\widehat{\lambda}_n$ estimate the quantities $G(\xi)$ and $\lambda$ from the observation of the $n$ first loss events. We establish that the distance between $\mathcal{R}$ and its estimate $\widehat{\mathcal{R}}_n$ is directly related to the estimation error of $\widehat{G}_n$ and $\widehat{\lambda}_n$.

\begin{proposition}\label{prop:maj:R}
Under assumption $(C^\lambda_1)$, the following inequality {almost surely} holds,
$$\sup_{(x,y)\in [1,\infty)\times[0,\infty)}\left|\mathcal{R}(x,y)-\widehat{\mathcal{R}}_n(x,y)\right|
\leq \left\|G-\widehat{G}_n\right\|_\infty + \frac{1}{\lambda_\ast}\left(4e^{-1}\frac{\lambda^\ast }{\lambda_\ast}+1\right)\left\|\widehat{G}_n\right\|_\infty \left|\lambda-\widehat{\lambda}_n\right|.$$
%$\prob$-a.s.
\end{proposition}

\noindent
This result has the following corollary.
\begin{corollary}
\label{prop:estimation:R}
Under assumptions $(C^\lambda_{1,2})$ and $(C^G_1)$,
the estimator $\widehat{\mathcal{R}}_n$ converges towards $\mathcal{R}$ in probability uniformly in $(x,y)\in[1,\infty)\times[0,\infty)$,
$$
\forall\,\e>0,\quad\lim_{n\to\infty}\prob\left(\sup_{(x,y)\in [1,\infty)\times[0,\infty)}\left|\mathcal{R}(x,y)-\widehat{\mathcal{R}}_n(x,y)\right|\geq\e\right)=0.
$$
In addition, the rate of convergence may be obtained from $(\ref{eq:rate:convergence})$.
\end{corollary}

\subsection{Absorption probability}\label{subs:abs}

The previous results on the estimation of the transition density $\mathcal{R}$ allow us to estimate the absorption probability of the PDMP $(X_t)_{t\geq0}$. When the trajectory starts from $X_0=x>1$, this probability is defined by
$$p(x) = \prob(X_t \in\Gamma\,\text{for some $t$}\,|\,X_0=x),$$
where $\Gamma=[0,1]$ is an absorbing set. We state that $p(x)$ may be found as a solution of an integral equation.
% As in \cite{TJRI}

\begin{proposition}\label{prop:eq:int}
$p(x)$ is a solution of the following integral equation,
\begin{equation}\label{eq:int:r}
p(x) = \int_0^1 \mathcal{R}(x,y)dy + \int_1^{+\infty} p(y)\mathcal{R}(x,y)dy.
\end{equation}
{The uniqueness of the solution to $(\ref{eq:int:r})$ will be discussed in Theorem \ref{prop:estim:r}.}
\end{proposition}
\proof First, we propose to rewrite $p(x)$ from the Markov chain $(Z_n)_{n\geq0}$,
$$p(x) = \prob(Z_n\in \Gamma\,\text{for some $n$}\,|\,Z_0=x) .$$
In addition, we have
$$p(x) = \prob(Z_1\in \Gamma\,|\,Z_0=x)\,+\,\prob(Z_1\notin\Gamma,\,Z_n\in\Gamma\,\text{for some $n\geq2$}\,|\,Z_0=x).$$
Together with
$$\prob(Z_1\notin\Gamma,\,Z_n\in\Gamma\,\text{for some $n\geq2$}\,|\,Z_0=x) = \int_{1}^\infty \prob(Z_n\in\Gamma\,\text{for some $n\geq2$}\,|\,Z_1=y)\,\mathcal{R}(x,y)\,dy,$$
{and the Markov property for $(Z_n)_{n\geq1}$}, this shows $(\ref{eq:int:r})$.\hfill$\Box$

%\noindent \textcolor{green}{The uniqueness of a solution to (\ref{eq:int:r}) will be discussed in Theorem \ref{prop:estim:r}.} 

\noindent
By virtue of $(\ref{eq:int:r})$, we propose to estimate $p(x)$ by the unique solution of the estimated integral equation
\begin{equation}\label{eq:int:r:hat}
\widehat{p}_n(x) = \int_0^1 \widehat{\mathcal{R}}_n(x,y)dy + \int_1^{+\infty} \widehat{p}_n(y)\widehat{\mathcal{R}}_n(x,y)dy,
\end{equation}
which satisfies both conditions
$$
\lim_{x\searrow 1}\widehat{p}_n(x)=1\qquad\text{and}\qquad\lim_{x\to\infty}\widehat{p}_n(x)=0.
$$
Nevertheless, the above equation is not in a proper form to compute $\widehat{p}_n$. As a consequence, we propose to solve numerically this estimated equation. On the space $L^1(1,\infty)$, endowed with its usual norm denoted by $\|\cdot\|$, we define the operator 
\begin{equation}\label{eq:operator:Kn}
\widehat{K}_n\,:\,h\mapsto \int_1^{+\infty} h(y)\,\widehat{\mathcal{R}}_n(x,y)\,dy,
\end{equation}
and we introduce the following additional notation,
\begin{equation}\label{eq:sn}
{\widehat{s}_n(\cdot)} = \int_0^1 \widehat{\mathcal{R}}_n(\cdot,y)dy .
\end{equation}
Thus, the equation (\ref{eq:int:r:hat}) may be rewritten as a Fredholm equation of the second kind on the space $L^1(1,\infty)$ \citep{NIETHAMMER},
\begin{equation}\label{eq:int:r:hat:fred}
\widehat{p}_n-\widehat{K}_n\widehat{p}_n=\widehat{s}_{n}.
\end{equation} 
As is well known, one may approximate a solution of (\ref{eq:int:r:hat:fred}) by the quantity
\begin{equation}\label{eq:pnm}
\widehat{p}_{n,m}=\sum_{k=0}^m \widehat{K}^k_n\widehat{s}_n,
\end{equation}
as long as $\|\widehat{K}_n\|<1$, this condition being ensured by an additional condition, as stated in the following theorem of convergence of $\widehat{p}_{n,m}$ towards $p$.
\begin{theorem}\label{prop:estim:r}
Under the conditions ($C^\lambda_{1,2}$) and ($C^G_2$), and the additional assumption  $\int_0^1G(u)\,u^{-1}du<1+r/\lambda$, the equation (\ref{eq:int:r}) has a unique solution and moreover, $\|\widehat{p}_{n,m}-p\|$ tends to $0$ in probability when $n$ and $m$ go to infinity.
%for all $\e>0$, there exists $N,M\in\N$ such that for all $n\geq N$ and $m\geq M$, there is some ``high-probability'' set $\Omega_n\subset \Omega$ (depending only on $n$) with $\prob(\Omega\setminus\Omega_n)\leq\e$ such that
%$$
%\|\widehat{p}_{n,m}-p\|\leq\e
%$$  
%on $\Omega_n$.
\end{theorem}
%Remark that the mode of convergence in Theorem \ref{prop:estim:r} is slightly weaker than convergence in probability. However, it still implies convergence in law.

\subsection{Distribution of the hitting time}
\label{subs:hit}

\noindent We now proceed to the estimation of $t_m(x)$, the probability for the process $(X_t)_{t\geq0}$ starting from $X_0=x$ to be absorbed at jump $m$. For $x>1$ we have $t_1(x)=\prob(Z_1\in\Gamma\,|\,Z_0=x)$ {and, as $\Gamma$ is absorbing}, for $m\geq2$,
$$
t_m(x)=\prob(Z_m\in\Gamma,\, {Z_{m-1}\notin\Gamma}\,|\,Z_0=x).
$$
We state in the following result that this sequence satisfies a recurrence relation.

\begin{proposition}\label{prop_t_m_exact}
For $x>1$, the functional sequence $(t_m)_{m\geq1}$ satisfies $t_1(x)=\int_0^1 \mathcal{R}(x,y)\, dy$ and the recursion relation,
\begin{equation}
\forall\,m\geq2,~t_m(x)=\int_1^\infty t_{m-1}(y)\mathcal{R}(x,y)\, dy.
\end{equation}
\end{proposition}

\proof
The proof follows the same reasoning as in the proof of Proposition \ref{prop:eq:int}.
\hfill$\Box$

\noindent
Following the same approach as in Subsection \ref{subs:abs}, we propose to estimate the functional sequence $(t_m)_{m\geq1}$ by the recursive procedure,
$$\widehat{t}_{1,n}(x)=\int_0^1 \widehat{\mathcal{R}}_n(x,y)\, dy,$$
and for $m\geq2$,
$$
\widehat{t}_{m,n}(x)=\int_1^\infty \widehat{t}_{m-1,n}(y)\,\widehat{\mathcal{R}}_n(x,y)\, dy.
$$
Using the operator $\widehat{K}_n$ defined in $(\ref{eq:operator:Kn})$ and the notation $(\ref{eq:sn})$, this recursion relation is closed to give
\begin{equation}\label{eq:tnm}
\widehat{t}_{m,n}=\widehat{K}^{m-1}_n\widehat{s}_n.
\end{equation}

\begin{theorem}\label{prop:estim:t}
Under the conditions ($C^\lambda_{1,2}$) and ($C^G_2$), and the additional assumption $\int_0^1G(u)\,u^{-1}du<1+r/\lambda$, then, for any integer $m$, $\|\widehat{t}_{m,n}-t_m\|$ tends to $0$ in probability when $n$ goes to infinity.

%for all $\e>0$ and $m\in\N$, there exists $N\in\N$ such that for all $n\geq N$, there is some ``high-probability'' set $\Omega_n\subset \Omega$ with $\prob(\Omega\setminus\Omega_n)\leq\e$ such that
%$$
%\|\widehat{t}_{m,n}-t_m\|\leq\e
%$$  
%on $\Omega_n$.
\end{theorem}

\noindent
We give the relation between the functional sequence $(\widehat{t}_{m,n})_{m\geq1}$ and the estimate $\widehat{p}_{n,m}$ of the absorption probability in the following remark.

\begin{remarque}
The estimation procedures for $p$ and $(t_m)_{m\geq1}$ may be carried out at the same time. In light of $(\ref{eq:sn})$, $(\ref{eq:pnm})$ and $(\ref{eq:tnm})$, we have
$$ \widehat{p}_{n,m} = \widehat{s}_n + \sum_{k=1}^m \widehat{t}_{k,n} .$$
As a consequence, the estimation of the absorption probability $p$ from the estimated sequence $(\widehat{t}_{m,n})_{m\geq1}$ does not require extra calculations.
\end{remarque}

%\subsection{Discussion about the hypotheses}
%
%In this section, we show that Assumptions \ref{hyp:G}, ($C^\lambda_1$), ($C^\lambda_2$) and ($C^G_3$) are satisfies by the ?? and Parzen-Rosenblatt estimators.

\section{Numerical illustration}
\label{sec:simulations}

This part of the paper is dedicated to some numerical illustrations of our main convergence results stated in the previous section. All the simulations have been implemented in the \verb+R+ language, which is commonly used in the statistical community, with an extensive use of the \verb+integrate+ function (numerical integration routine with adaptive quadrature of functions). As an example in our simulations, we choose for the probability density function $G$ the following power function, $G(u)=11u^{10}$ for any $u\in[0,1]$. This density function charges the interval $[0.8,1]$ at more than $90\%$. This means that the process is weakly affected by a fragmentation event. For the jump rate we choose $\lambda=1$ and for the growth rate $r=1$. Then,
$$
%\|K\| \leq
\frac{\lambda}{\lambda+r} \int_0^1G(u)\,u^{-1}\,du\,=\, 0.55\,<\,1,
$$
so that we are in the scope of application of Theorems \ref{prop:estim:r} and \ref{prop:estim:t}. We propose to illustrate our theoretical results Corollary \ref{prop:estimation:R} and Theorems \ref{prop:estim:r} and \ref{prop:estim:t} from the observation of different numbers of data ($n=50,\,75$ and $100$). In addition, we always present the distribution of our estimates from a fixed number of data over $100$ replicates of the numerical experiment.

For these simulation experiments, we choose to estimate the density $G(x)$ by the Parzen-Rosenblatt estimator $\widehat{G}_n^{PR}(x)$ defined by
\begin{equation*}
\widehat{G}_n^{PR}(x) = \frac{1}{n h_n} \sum_{i=1}^n \mathbb{K}\left(\frac{Y_i-x}{h_n}\right),
\end{equation*}
where $\mathbb{K}$ is the Gaussian kernel and the parameter $h_n$ is the bandwidth. The estimator is computed from the \verb+R+ function \verb+density+ with an optimal choice of the bandwidth parameter. In addition, $\lambda$ is estimated from the observations $S_i$'s by the truncated maximum likelihood estimator $\widehat{\lambda}_n^{ml}$ defined in $(\ref{eq:tmle})$. These estimates satisfy the conditions that we impose in the paper.

First, we present some simulation results for the transition kernel $\mathcal{R}(x,y)$ {(see Figures \ref{fig:Rx2}, \ref{fig:R2y} and \ref{fig:R:ISE})}. The transition kernel is not really a quantity of interest in the model in contrary to the rate and measure of jumps $\lambda$ and $G(u)du$. Nevertheless, the kernel appears when we want to compute the probability of hitting, or the hitting time of, $\Gamma$. This is therefore required to be able to estimate $\mathcal{R}(x,y)$ in our approach. Recall the definition $(\ref{eq:estim:R})$ of the estimator of $\mathcal{R}$ from $\widehat{\lambda}_n$ and $\widehat{G}_n^{PR}$. In Figure \ref{fig:Rx2} are displayed {the trajectory of $\mathcal{R}(x,2)$ and its estimates for $1\leq x\leq 4$ from $n=50$, $75$ and $100$ data as well as the pointwise error with boxplots over $100$ replications between $\mathcal{R}(\cdot,2)$ and $\widehat{\mathcal{R}}_{100}(\cdot,2)$ within the interval $[1,4]$. Figure \ref{fig:R2y} presents the same numerical results for the estimation of $\mathcal{R}(2,y)$, $1\leq y\leq4$. Notice that, according to Figures \ref{fig:Rx2} and \ref{fig:R2y}, the pointwise error in the estimation of $\mathcal{R}(\cdot,2)$ and $\mathcal{R}(2,\cdot)$ is maximum around $2$. This may be explained by the presence of a singularity in $2$ for both these functions.} The corresponding integrated square errors are given in Figure \ref{fig:R:ISE}. In both cases, we observe a decrease in the error when the number of data grows, despite the low number of data. However, this is not very surprising here since the transition kernel is estimated from its exact expression (see Proposition \ref{prop:def:ker}), substituting $\lambda$ by $\widehat{\lambda}_n$ and $G$ by $\widehat{G}_n^{PR}$.

\begin{center}[\,Figures \ref{fig:Rx2}, \ref{fig:R2y} and \ref{fig:R:ISE} here.\,]\end{center}

Now, we proceed to the simulation of the estimation of $p(x)$, the probability for the process $(X_t)_{t\geq0}$ to be absorbed by $\Gamma=[0,1]$ starting form $x>1$. This is, with the time of absorption, one of the two main quantities of interest in the model. Indeed, for the ruin theoretical model of Section \ref{sec_ex_eco}, $p(x)$ corresponds to the probability to be {ruined} starting from some capital $x$. For the Malthusian evolution of Section \ref{sec_ex_bio}, $p(x)$ is the probability for a population of initial size $x$ to extinct. Nevertheless, we can not compute directly the function of interest $p$. As a consequence, we propose to compare $\widehat{p}_{n,m}$ and the numerical approximation $p_m=\sum_{k=0}^mK^ks$ of $p$, where the operator $K$ is defined in $(\ref{eq:operatorK})$ and $s(\cdot)=\int_0^1\mathcal{R}(\cdot,y)\,dy$. Roughly speaking, $K$ and $s$ are the deterministic limits of the estimates $\widehat{K}_n$ and $\widehat{s}_n$ presented in $(\ref{eq:operator:Kn})$ and $(\ref{eq:sn})$. The error in $L^1(1,\infty)$-norm between $p$ and $p_m$ satisfies
$$
\|p-p_m\|\leq \|s\|\frac{\|K\|^{m+1}}{1-\|K\|}.
$$
Together with the chosen numerical values and $m=10$, we have $\|p-p_m\|\leq 1.6\times10^{-4}$. Consequently, the numerical error due to the approximation of $p$ is very low and does not affect our comparison results presented in the sequel.

Recall that our approximation of $p$ is given by $\widehat{p}_{n,m}=\sum_{k=0}^m \widehat{K}^k_n\widehat{s}_n$. In the simulations, we compare $\widehat{p}_{n,m}$ with $p_m$ for $m=10$ and $n=50$, $70$ and $100$ data. Figure \ref{fig:p:curve} displays the shape of $p_m$ and $\widehat{p}_{n,m}$ as well as {the boxplots of the punctual error between the curves $p_m$ and $p_{100,m}$}. The corresponding integrated square error is presented in Figure \ref{fig:p:ise}. A decrease in the error is observed when $n$ grows. Note that the error is already small for $n=50$ and seems to behave quite well despite the successive application of the kernel $\widehat{K}_n$.

\begin{center}[\,Figures \ref{fig:p:curve} and \ref{fig:p:ise} here.\,]\end{center}

Finally, we go on with the estimation of $t_m(x)$, the probability for the process $(X_t)_{t\geq0}$ starting from $x$, to be absorbed at jump $m$. The quantity $t_m(x)$ is an important feature of the model and provides additional information to that given by $p(x)$. Remark that according to Proposition \ref{prop_t_m_exact}, $t_m$ may be computed in an exact way {contrary to} $p(x)$. There is therefore no numerical error in this case (if we do not consider the numerical errors introduced by the computation of the kernel integrals). Thus, we compare directly $t_m(x)$ with its estimator $\widehat{t}_{m,n}(x)$ given by equation (\ref{eq:tnm}). At first, we notice that the estimation of the probability of absorption $\widehat{p}_{n,m}$ and the estimation of the times at which an absorption occurs $\widehat{t}_{m,n}(x)$ are related through the formula,
$$
\widehat{p}_{n,m}=\sum_{k=0}^{m}\widehat{t}_{k+1,n}.
$$
Therefore, in our previous computations of $\widehat{p}_{n,m}$, we already have computed the quantities $\widehat{t}_{m,n}$ and no further calculations are required. In Figure \ref{fig:t:ise}, we present the integrated square error between $t_m$ and its estimate $\widehat{t}_{m,n}$ from the observation of $n=50,\,75$ or $100$ random loss events and for $m=1,\,2,\,3$ and $4$, that is for the four first absorption times. There is a decrease of the error when $n$ grows for each value of $m$. Quantitatively, this does not make sense to compare the error for $m=2$ and $m=4$ since, as displayed in Figure \ref{fig:t}, the order of magnitude of the estimated probabilities is not at all the same. Figure \ref{fig:t} presents the distribution of the hitting time of $\Gamma$, $t_m(x)$, for $x=1.1$ and $m=1,\,\dots,\,6$, and also the distribution of its estimates $\widehat{t}_{m,n}$ from the observation of $n=50,\,75$ or $100$ random loss events. More precisely, in this figure is represented the mean of the estimators together with the first and third quartiles, over $100$ replications. Once again, a decrease in the error was observed when $n$ grows showing that the law of the hitting times of $\Gamma$ is well estimated. These results, coupled with the estimate of $p(x)
$, give all the interesting information in the study of this model. In all the procedure, the estimates are of high quality despite the low number of data used, in particular with respect to the sample sizes used in \citep{AzaisESAIM14,AzaisSJS14,Doumic11,Doumic12}.

\begin{center}[\,Figures \ref{fig:t:ise} and \ref{fig:t} here.\,]\end{center}

\noindent\textbf{Acknowledgments:}~~The referees deserve thanks for careful reading of the original version of the manu\-script and many helpful suggestions for improvement in the article. The authors also acknowledge Alexandre Boumezoued for fruitful discussions about hybrid processes and Poisson random measures.

\appendix
\section{Some technical lemmas}\label{sec:annexe}

This part is dedicated to the presentation of some technical results which will be useful in the proofs of our main results presented in Appendix \ref{app:proofs}. For convenience, we use in the sequel the following {notation}. For $\lambda>0$, $x\geq1$, $y\geq0$ {and $u\in[0,1]$}, we define
\begin{equation}\label{nota:proof}
\alpha_\lambda(x) = \frac{(x-1)^{\lambda/r}}{r},\quad\beta_\lambda(y,u)=u^{\lambda/r} (y-u)^{-\lambda/r-1}\quad\text{and}\quad f_\lambda(x,y)=\alpha_\lambda(x)\int_{0}^{y/x\wedge 1}\beta_\lambda(y,u)du.
\end{equation}

\begin{lemma}\label{lem:Yn:iid}
{The sequence $(Y_n)_{n\geq1}$ has $G$ as common distribution, is independent and independent of the interarrival times $(S_n)_{n\geq1}$.}
\end{lemma}
\proof
{For any integer $n$, the $\sigma$-algebra $\sigma(X_0,\,Y_1,\,\dots,\,Y_{n-1},\,S_1,\,\dots,\,S_n)$ is denoted $\mathcal{F}_{n-1}$. First, let us notice that the post-jump location $Z_{n-1}$ is $\mathcal{F}_{n-1}$-measurable. By the expression of the transition kernel $\mathcal{Q}$ $(\ref{eq:def:Q})$, for any measurable function $\varphi$, we have
\begin{eqnarray*}
\esp\left[\varphi(Y_n)\,\big|\,\mathcal{F}_{n-1}\right] &=&\esp\left[\varphi\left(\frac{Z_n}{\Phi(Z_{n-1},S_n)}\right)\,\Bigg|\,\mathcal{F}_{n-1}\right]\\
~&=& \int \varphi\left(\frac{\zeta}{\Phi(Z_{n-1},S_n)}\right)\,{\Phi(Z_{n-1},S_n)}^{-1}\,G\left(\frac{\zeta}{\Phi(Z_{n-1},S_n)}\right) d\zeta\\
~&=& \int \varphi(y)G(y) dy,
\end{eqnarray*}
by the change of variables $\zeta=y\Phi(Z_{n-1},S_n)$. This yields the expected result.}\hfill$\Box$

\begin{lemma}\label{lem:inverseflow}
For any $x>1$, $y\geq x$, the deterministic flow $(\ref{eq:def:flow})$ satisfies ${\Phi(x,t)}=y$ if and only if $t=\frac{1}{r}\log\left(\frac{y-1}{x-1}\right)$.
\end{lemma}
\proof This result is obvious.\hfill$\Box$
\begin{lemma}\label{lem:lip}
%Recall the notations (\ref{nota:proof}).
For any ${\lambda},\lambda_1,\lambda_2\in[\lambda_\ast,\lambda^\ast]$, $x\geq1$, $y\geq0$ and $u\in(0,y/x)$, we have
\begin{align*}
&(i)~f_\lambda(x,y)\,\leq\, \frac{1}{\lambda}\mathbb{1}_{\{y<x\}}\left[1-\left(\frac{x-1}{x}\right)^{\lambda/r}\right]\,+\,\frac{1}{\lambda}\mathbb{1}_{\{y\geq x\}}\left[\left(\frac{x-1}{y-1}\right)^{\lambda/r}-\left(\frac{x-1}{y}\right)^{\lambda/r}\right]\,\leq\, \frac{1}{\lambda},\\
&(ii)~\Big|\alpha_{\lambda_1}(x)\beta_{\lambda_1}(y,u)-\alpha_{\lambda_2}(x)\beta_{\lambda_2}(y,u)\Big|\,\leq\,\frac{1}{r^2}\frac{1}{y-u}\left(\frac{u(x-1)}{y-u}\right)^{\lambda_\ast/r}\log\left(\frac{u(x-1)}{y-u}\right)\big|\lambda_1-\lambda_2\big|,\\
&(iii)~ \sup_{x\geq1,y\geq0}\int_0^{y/x\wedge 1}\Big|\alpha_{\lambda_1}(x)\beta_{\lambda_1}(y,u)-\alpha_{\lambda_2}(x)\beta_{\lambda_2}(y,u)\Big|du\,\leq\, \frac{4e^{-1}}{\lambda^2_\ast}\big|\lambda_1-\lambda_2\big|.
\end{align*}
\end{lemma}
\proof
We begin with the inequality (i). Let $(x,y)$ be in $[1,\infty)\times[0,\infty)$. For any $u\in[0,y/x\wedge 1]$ we have
\[
0\leq u^{\lambda/r}(y-u)^{-\lambda/r-1}\,\leq\,  \left(y/x\wedge 1\right)^{\lambda/r}(y-u)^{-\lambda/r-1}.
\]
Therefore,
\[
f_\lambda(x,y)\leq\frac{1}{r} \left(y/x\wedge 1\right)^{\lambda/r}(x-1)^{\lambda/r} \int_{0}^{y/x\wedge 1} (y-u)^{-\lambda/r-1}du.
\]
Computing the integral leads to
\[
f_\lambda(x,y)\leq \frac{1}{\lambda}\left(y/x\wedge 1\right)^{\lambda/r}(x-1)^{\lambda/r}\left[\left(y-y/x\wedge 1\right)^{-\lambda/r}-y^{-\lambda/r}\right].
\]
We now split the latter term in two using the elementary fact that $1=\mathbb{1}_{\{y<x\}}+\mathbb{1}_{\{y\geq x\}}$. This yields
\begin{align*}
f_{\lambda}(x,y)\,\leq\,&~\frac{1}{\lambda}\left(y/x\right)^{\lambda/r}(x-1)^{\lambda/r}\mathbb{1}_{\{y<x\}}\left[\left(y-y/x\right)^{-\lambda/r}-y^{-\lambda/r}\right]\\
&+\frac{1}{\lambda} (x-1)^{\lambda/r}\mathbb{1}_{\{y\geq x\}}\left[\left(y- 1\right)^{-\lambda/r}-y^{-\lambda/r}\right]\\
\,=\,&~\frac{1}{\lambda}\mathbb{1}_{\{y<x\}}\left[1-\left(\frac{x-1}{x}\right)^{\lambda/r}\right]\\
&+\frac{1}{\lambda}\mathbb{1}_{\{y\geq x\}}\left[\left(\frac{x-1}{y-1}\right)^{\lambda/r}-\left(\frac{x-1}{y}\right)^{\lambda/r}\right].
\end{align*}
Using the fact that when $y\geq x$, $\frac{x-1}{y-1}\leq 1$ and noticing that the two terms $\frac{x-1}{y}$ and $\frac{x-1}{x}$ are non negative, we obtain:
\[
f_{\lambda}(x,y)\leq\frac{1}{\lambda}\mathbb{1}_{\{y< x\}}+\frac{1}{\lambda}\mathbb{1}_{\{y\geq x\}}=\frac{1}{\lambda}.
\]
We go on with the second part of the lemma. One may derivate with respect to $\lambda$ to obtain
$$
\partial_\lambda \alpha_\lambda(x)\beta_\lambda(y,u)=\frac{1}{r^2}\frac{1}{y-u}\left(\frac{u(x-1)}{y-u}\right)^{\lambda/r}\log\left(\frac{u(x-1)}{y-u}\right).
$$
Then, we use that for $x\geq1$, $y\geq0$ and $u\in[0,y/x\wedge 1]$, one has $X=\frac{u(x-1)}{y-u}\in(0,1]$ such that $X^{\lambda/ 2r}|\log X|\leq \frac{2r}{\lambda}e^{-1}$. This fact yields
$$
\big|\partial_\lambda \alpha_\lambda(x)\beta_\lambda(y,u)\big|\,\leq\,\frac{2e^{-1}}{r\lambda}\frac{1}{y-u}\left(\frac{u(x-1)}{y-u}\right)^{\lambda/2r}\,\leq\,\frac{2e^{-1}}{r\lambda_\ast}\frac{1}{y-u}\left(\frac{u(x-1)}{y-u}\right)^{\lambda_\ast/2r}.
$$
Notice that the last inequality is uniform in $\lambda\in[\lambda_\ast,\lambda^\ast]$. This proves the second assertion (ii). For the third one, using the mean value theorem and similar calculations as above, we obtain that
\begin{align*}
\int_0^{y/x\wedge 1}\Big|\alpha_{\lambda_1}(x)\beta_{\lambda_1}(y,u)-\alpha_{\lambda_2}(x)\beta_{\lambda_2}(y,u)\Big|du\leq&~\mathbb{1}_{\{y<x\}}\frac{4e^{-1}}{\lambda^2_\ast}\left[1-\left(\frac{x-1}{x}\right)^{\lambda_\ast/2r}\right]\big|\lambda_1-\lambda_2\big|\\
&+\mathbb{1}_{\{y\geq x\}}\frac{4e^{-1}}{\lambda^2_\ast}\left[\left(\frac{x-1}{y-1}\right)^{\lambda_\ast/2r}-\left(\frac{x-1}{y}\right)^{\lambda_\ast/2r}\right]\big|\lambda_1-\lambda_2\big|.
\end{align*}
Using the fact that when $y\geq x$, $\frac{x-1}{y-1}\leq 1$ and noticing that the two terms $\frac{x-1}{y}$ and $\frac{x-1}{x}$ are non negative, we obtain
\begin{align*}
\int_0^{y/x\wedge 1}\Big|\alpha_{\lambda_1}(x)\beta_{\lambda_1}(y,u)-\alpha_{\lambda_2}(x)\beta_{\lambda_2}(y,u)\Big|du~~\leq&~~\mathbb{1}_{\{y<x\}}\frac{4e^{-1}}{\lambda^2_\ast}\big|\lambda_1-\lambda_2\big|+\mathbb{1}_{\{y\geq x\}}\frac{4e^{-1}}{\lambda^2_\ast}\big|\lambda_1-\lambda_2\big|\\
=&~~\frac{4e^{-1}}{\lambda^2_\ast}\big|\lambda_1-\lambda_2\big|.
\end{align*}
The result follows.
\hfill$\Box$
%

%%
%\begin{lemma}\label{lemma:calR:G}
%Under Assumption \ref{hyp:G}, the transition kernel $\mathcal{R}$ satisfies
%\[
%\sup_{[1,\infty)\times[0,\infty)}\mathcal{R}\leq \sup_{[0,1]} G.
%\]
%\end{lemma}
%\proof
%Let $(x,y)$ be in $[1,\infty)\times[0,\infty)$. The conditional density $\mathcal{R}(x,y)$ given by $(\ref{eq:calcul:R})$ satisfies
%%Recall the definition of $\mathcal{R}$ when $x\geq1$,
%%\[
%%\mathcal{R}(x,y)=\frac{\lambda}{r}(x-1)^{\lambda/r} \int_{0}^{y/x\wedge 1} G(u) u^{\lambda/r}(y-u)^{-\lambda/r-1}du.
%%\]
%%Therefore,
%$$
%\mathcal{R}(x,y)\leq \lambda\sup_{u\in[0,1]} G(u) f_{\lambda}(x,y).
%$$
%Using the previous lemma, we get, uniformly in $x\geq1$ and $y\geq0$,
%\[
%\mathcal{R}(x,y)\,\leq\,\lambda\sup_{u\in[0,1]} G(u)\frac{1}{\lambda}\,=\,\sup_{u\in[0,1]} G(u).
%\]
%This shows the expected result.\hfill$\Box$

\begin{lemma}\label{lem-cond-K}
The following equality holds,
$$
\sup_{y\in[0,\infty)}\int_1^\infty \mathcal{R}(x,y)dx\,=\,\frac{\lambda}{\lambda+r}\int_0^1 \frac{G(u)}{u}\,du.
$$
\end{lemma}

\proof
By definition of $\mathcal{R}$ one may write
$$
\int_1^\infty\mathcal{R}(x,y)\,dx \, = \, \lambda\int_1^\infty \alpha_\lambda(x) \int_0^{y/x\wedge1}\beta_\lambda(y,u)\,G(u)\,du\,dx.
$$
In the above term, one may change the order of integration to integrate in $x$ at first. We obtain
\begin{eqnarray*}
\lambda \int_1^\infty\alpha_\lambda(x) \int_0^{y/x\wedge1}\beta_\lambda(y,u)G(u)\,du\,dx &=&\lambda\int_0^{1} \int_1^{y/u}\alpha_\lambda(x) \,dx\beta_\lambda(y,u)G(u)\mathbb{1}_{\{u\leq y\}}\,du\\
&=&\frac{\lambda}{r}\int_0^{1} \frac{1}{\lambda/r+1}\left(\frac yu -1\right)^{\lambda/r+1}\beta_\lambda(y,u)G(u)\mathbb{1}_{\{u\leq y\}}\,du\\
&=&\frac{\lambda}{\lambda+r}\int_0^{1}G(u)\mathbb{1}_{\{u\leq y\}}\,\frac{du}{u}.
\end{eqnarray*}
%\begin{align*}
%&\lambda \int_1^\infty\alpha_\lambda(x) \int_0^{y/x\wedge1}\beta_\lambda(y,u)G(u)\,du\,dx\\
%&=\lambda\int_0^{1} \int_1^{y/u}\alpha_\lambda(x) \,dx\beta_\lambda(y,u)G(u)1_{u\leq y}\,du\\
%&=\frac{\lambda}{r}\int_0^{1} \frac{1}{\lambda/r+1}\left(\frac yu -1\right)^{\lambda/r+1}\beta_\lambda(y,u)G(u)1_{u\leq y}\,du\\
%&=\frac{\lambda}{\lambda+r}\int_0^{1}G(u)1_{u\leq y}\,\frac{du}{u}.
%\end{align*}
The result follows.
\hfill$\Box$

\begin{lemma}\label{hatR-specialnorm}
Under Assumption $(C^\lambda_1)$, almost-surely the following inequality holds
$$
\sup_{y\in[0,\infty)}\int_1^\infty\!\left|\mathcal{R}(x,y)-\widehat{\mathcal{R}}_n(x,y)\right|\!dx\leq\frac{\lambda}{\lambda+r}\int_0^1\left|G(u)-\widehat{G}_n(u)\right|\frac{du}{u}+\lambda^\ast\!\!\left[4e^{-1}\frac{1}{\lambda^2_\ast }+\frac{1}{\lambda_\ast+r}\right]\!\int_0^1|\widehat{G}_n(u)|\frac{du}{u}\big|\lambda-\widehat{\lambda}_n\big|.
$$
\end{lemma}

\proof
One may write
\begin{eqnarray*}
\mathcal{R}(x,y) - \widehat{\mathcal{R}}_n(x,y) &=& \lambda \alpha_\lambda(x) \int_0^{y/x\wedge1}\beta_\lambda(y,u)\,\left[G(u)-\widehat{G}_n(u)\right]\,du\\
&&\!+\,\,\lambda \int_0^{y/x\wedge1} \widehat{G}_n(u)\,\left[\alpha_\lambda(x)\beta_\lambda(y,u)-\alpha_{\widehat{\lambda}_n}(x)\beta_{\widehat{\lambda}_n}(y,u)\right]\,du\\
&&\!+\,\,\left(\lambda-\widehat{\lambda}_n\right)\,\alpha_{\widehat{\lambda}_n}(x)\int_0^{y/x\wedge1}\beta_{\widehat{\lambda}_n}(y,u)\widehat{G}_n(u)du.
\end{eqnarray*}
Thus, for any $y\geq0$, we have
\begin{align*}
\int_1^\infty\left|\mathcal{R}(x,y) - \widehat{\mathcal{R}}_n(x,y)\right|\,dx&\,\,\leq\,\,\lambda \int_1^\infty\alpha_\lambda(x) \int_0^{y/x\wedge1}\beta_\lambda(y,u)\,\left|G(u)-\widehat{G}_n(u)\right|\,du\,dx\\
&\qquad+\,\lambda \int_1^\infty\int_0^{y/x\wedge1} |\widehat{G}_n(u)|\,\left|\alpha_\lambda(x)\beta_\lambda(y,u)-\alpha_{\widehat{\lambda}_n}(x)\beta_{\widehat{\lambda}_n}(y,u)\right|\,du\,dx\\
&\qquad+\,\left|\lambda-\widehat{\lambda}_n\right|\,\int_1^\infty\alpha_{\widehat{\lambda}_n}(x)\int_0^{y/x\wedge1}\beta_{\widehat{\lambda}_n}(y,u)|\widehat{G}_n(u)|\,du\,dx.
\end{align*}
In the three above terms, one may change the order of integration to integrate in $x$ at first. Let us deal with these terms separately. For the first term, we have
\begin{align*}
\lambda \int_1^\infty\alpha_\lambda(x) \int_0^{y/x\wedge1}\beta_\lambda(y,u)\,\big|G(u)-&\widehat{G}_n(u)\big|\,du\,dx\\
&=\,\,\lambda\int_0^{1} \int_1^{y/u}\alpha_\lambda(x) \,dx\beta_\lambda(y,u)\,\big|G(u)-\widehat{G}_n(u)\big|\mathbb{1}_{\{u\leq y\}}\,du\\
&=\,\,\frac{\lambda}{r}\int_0^{1} \frac{1}{\lambda/r+1}\left(\frac yu -1\right)^{\lambda/r+1}\beta_\lambda(y,u)\,\big|G(u)-\widehat{G}_n(u)\big|\mathbb{1}_{\{u\leq y\}}\,du\\
&=\,\,\frac{\lambda}{\lambda+r}\int_0^{1}\big|G(u)-\widehat{G}_n(u)\big|\mathbb{1}_{\{u\leq y\}}\,\frac{du}{u}.
\end{align*}
%\begin{align*}
%&\lambda \int_1^\infty\alpha_\lambda(x) \int_0^{y/x\wedge1}\beta_\lambda(y,u)\,\left|G(u)-\widehat{G}_n(u)\right|\,du\,dx\\
%&=\lambda\int_0^{1} \int_1^{y/u}\alpha_\lambda(x) \,dx\beta_\lambda(y,u)\,\left|G(u)-\widehat{G}_n(u)\right|1_{u\leq y}\,du\\
%&=\frac{\lambda}{r}\int_0^{1} \frac{1}{\lambda/r+1}\left(\frac yu -1\right)^{\lambda/r+1}\beta_\lambda(y,u)\,\left|G(u)-\widehat{G}_n(u)\right|1_{u\leq y}\,du\\
%&=\frac{\lambda}{\lambda+r}\int_0^{1}\left|G(u)-\widehat{G}_n(u)\right|1_{u\leq y}\,\frac{du}{u}.
%\end{align*}
For the third term, a similar calculation gives
$$\big|\lambda-\widehat{\lambda}_n\big|\,\int_1^\infty\alpha_{\widehat{\lambda}_n}(x)\int_0^{y/x\wedge1}\beta_{\widehat{\lambda}_n}(y,u)|\widehat{G}_n(u)|\,du\,dx\,\,=\,\,\big|\lambda-\widehat{\lambda}_n\big|\frac{\widehat{\lambda}_n}{\widehat{\lambda}_n+r}\int_0^{1}\left|\widehat{G}_n(u)\right|\,\mathbb{1}_{\{u\leq y\}}\frac{du}{u}.
$$
The most intricate term is the second. Using Lemma \ref{lem:lip}, we have
\begin{align*}
\lambda \int_1^\infty\int_0^{y/x\wedge1}\widehat{G}_n(u)&\,\Big|\alpha_\lambda(x)\beta_\lambda(y,u)-\alpha_{\widehat{\lambda}_n}(x)\beta_{\widehat{\lambda}_n}(y,u)\Big|\,du\,dx\\
&\leq\,\,\,\lambda \big|\lambda-\widehat{\lambda}_n\big|\int_1^\infty\int_0^{y/x\wedge1} \frac{2e^{-1}}{r\lambda_\ast}\frac{1}{y-u}\left(\frac{u(x-1)}{y-u}\right)^{\lambda_\ast/2r}\big|\widehat{G}_n(u)\big|\,du\,dx\\
&=\,\,\,\frac{\lambda 2e^{-1}}{r\lambda_\ast}\big|\lambda-\widehat{\lambda}_n\big|\int_0^1\frac{1}{y-u}\left(\frac{u}{y-u}\right)^{\lambda_\ast/2r}\big|\widehat{G}_n(u)\big|\int_1^{y/u}(x-1)^{\lambda_\ast/2r}dx\,\mathbb{1}_{\{u\leq y\}}\,du\\
&=\,\,\,\frac{4\lambda^\ast e^{-1}}{\lambda^2_\ast }\big|\lambda-\widehat{\lambda}_n\big|\int_0^1\big|\widehat{G}_n(u)\big|\mathbb{1}_{\{u\leq y\}}\,\frac{du}{u}.
\end{align*}
The result follows by aggregation of the three above estimates.
\hfill$\Box$
\section{Proofs of the main results}
\label{app:proofs}

This section gathers the proofs of the different propositions stated in Section \ref{sec:main:res}.
%For convenience, we use in the sequel the following notations. For $\lambda>0$, $x\geq1$ and $y\geq0$, we define,
%\begin{equation}\label{nota:proof}
%\alpha_\lambda(x) = \frac{(x-1)^{\lambda/r}}{r},\,\beta_\lambda(y,u)=u^{\lambda/r} (y-u)^{-\lambda/r-1}\,\,\text{and}\,\, f_\lambda(x,y)=\alpha_\lambda(x)\int_{0}^{y/x\wedge 1}\beta_\lambda(y,u)du.
%\end{equation}

\subsection{Proof of Proposition \ref{prop:def:ker}}

In both cases $x>1$ and $x\leq1$, we have
\begin{equation}\label{eq:R}
\mathcal{R}(x,dy)=\left[\int_{\R_+}\frac{1}{z}G\left(\frac{y}{z}\right)\mathcal{S}(x,dz)\right]\,dy,
\end{equation}
where the conditional distribution $\mathcal{S}(x,dz)$ is defined from its cumulative version,
$$\mathcal{S}(x,(-\infty,z])=\prob\left({\Phi(Z_{n-1},S_n)}\leq z\,|\,Z_{n-1}=x\right). $$
For $x\leq 1$, from $(\ref{eq:def:flow})$, we have $\mathcal{S}(x,dz)=\delta_{\{x\}}(dz)$. This shows $(\ref{eq:calcul:R})$ for $x\leq1$.
If $x>1$, for any $z\geq x$,
$$\mathcal{S}(x,(-\infty,z]) = \prob\left(S_n\leq\frac{1}{r}\log\left(\frac{z-1}{x-1}\right) \right),$$
according to Lemma \ref{lem:inverseflow}. As a consequence, {we have}
\begin{equation}\label{eq:calcul:S}
\mathcal{S}(x,dz)=\frac{\lambda}{r} \,\frac{ (x-1)^{\lambda/r}}{(z-1)^{\lambda/r+1}}\,\mathbb{1}_{[x,\infty)}(z)\, dz.
\end{equation}
Since $G$ is a probability density function on $[0,1]$, together with $(\ref{eq:calcul:S})$, we may re-write $(\ref{eq:R})$ as
\begin{equation*}
\mathcal{R}(x,dy)=\frac{\lambda}{r} \,(x-1)^{\lambda/r}\,\left[\int_{x\vee y}^{+\infty} G\left(\frac{y}{z}\right) \frac{(z-1)^{-\lambda/r-1}}{z}dz\right]\,dy
\end{equation*}
By the change of variable $u=y/z$, we obtain
$$\mathcal{R}(x,y)=\frac{\lambda}{r} (x-1)^{\lambda/r}\left[\int_{0}^{y/x\wedge1}G(u) \left(\frac{y-u}{u}\right)^{-\lambda/r-1}\,u^{-1}\,du\right].$$
This shows the result $(\ref{eq:calcul:R})$ for $x>1$.

\subsection{Proof of Proposition \ref{prop:maj:R}}

Let $n\in\N$, $\lambda\in[\lambda_\ast,\lambda^\ast]$, $x\geq1$ and $y\geq0$. We work $\omega$ by $\omega$ so that the desired almost-sure inequality will follow. Recall that by equations (\ref{eq:calcul:R}) and (\ref{eq:estim:R}) together with the notations $(\ref{nota:proof})$,
$$\mathcal{R}(x,y)=\lambda \alpha_\lambda(x)\int_0^{y/x\wedge 1}\beta_\lambda(y,u)G(u)du$$
and
$$\widehat{\mathcal{R}}_n(x,y) = \widehat{\lambda}_n \alpha_{\widehat{\lambda}_n}\int_0^{y/x\wedge 1}\beta_{\widehat{\lambda}_n}(y,u)\widehat{G}_n(u)du .$$
By an elementary rearranging, one may write
\begin{align*}
\mathcal{R}(x,y) - \widehat{\mathcal{R}}_n(x,y) =&~ \lambda \alpha_\lambda(x) \int_0^{y/x\wedge1}\beta_\lambda(y,u)\,\left[G(u)-\widehat{G}_n(u)\right]\,du\\
&+\lambda \int_0^{y/x\wedge1} \widehat{G}_n(u)\,\left[\alpha_\lambda(x)\beta_\lambda(y,u)-\alpha_{\widehat{\lambda}_n}(x)\beta_{\widehat{\lambda}_n}(y,u)\right]\,du\\
&+ \left(\lambda-\widehat{\lambda}_n\right)\,\alpha_{\widehat{\lambda}_n}(x)\int_0^{y/x\wedge1}\beta_{\widehat{\lambda}_n}(y,u)\widehat{G}_n(u)du.
\end{align*}
We deal with the three above terms separately. For the first term we have
$$
\lambda \alpha_\lambda(x) \left|\int_0^{y/x\wedge1}\beta_\lambda(y,u)\,\left[G(u)-\widehat{G}_n(u)\right]\,du\right|\,\leq\,\lambda \alpha_\lambda(x) \int_0^{y/x\wedge1}\beta_\lambda(y,u)\,du\,\|G-\widehat{G}_n\|_\infty.
$$
Thus, by the first part of Lemma \ref{lem:lip}, {we obtain}
$$
\left|\lambda \alpha_\lambda(x) \int_0^{y/x\wedge1}\beta_\lambda(y,u)\,\left[G(u)-\widehat{G}_n(u)\right]\,du\right|\,\leq\,\lambda f_\lambda(x,y)\|G-\widehat{G}_n\|_\infty\,\leq\, \|G-\widehat{G}_n\|_\infty.
$$
Now for the second term, using this time the third part of Lemma \ref{lem:lip},
\begin{align*}
\lambda \Bigg|\int_0^{y/x\wedge1} \widehat{G}_n(u)\,\Big[\alpha_\lambda(x)\beta_\lambda(y,u)&-\alpha_{\widehat{\lambda}_n}(x)\beta_{\widehat{\lambda}_n}(y,u)\Big]\,du\Bigg|\\
&\leq~\lambda\|\widehat{G}_n(u)\|_\infty \int_0^{y/x\wedge1} \left|\alpha_\lambda(x)\beta_\lambda(y,u)-\alpha_{\widehat{\lambda}_n}(x)\beta_{\widehat{\lambda}_n}(y,u)\right|\,du\\
&\leq~\lambda\|\widehat{G}_n(u)\|_\infty\frac{4e^{-1}}{\lambda^2_\ast}|\lambda-\widehat{\lambda}_n|\\
&\leq~\frac{4\lambda^\ast e^{-1}}{\lambda^2_\ast}\|\widehat{G}_n(u)\|_\infty|\lambda-\widehat{\lambda}_n|.
\end{align*}
For the last term we have, using again the first part of Lemma \ref{lem:lip},
$$
\left|\lambda-\widehat{\lambda}_n\right|\,\alpha_{\widehat{\lambda}_n}(x)\int_0^{y/x\wedge1}\beta_{\widehat{\lambda}_n}(y,u)\widehat{G}_n(u)du\,\leq\,\left|\lambda-\widehat{\lambda}_n\right|\,\|\widehat{G}_n\|_\infty f_{\widehat{\lambda}_n}(x,y)\,\leq\,\left|\lambda-\widehat{\lambda}_n\right|\,\frac{\|\widehat{G}_n\|_\infty}{\lambda_\ast}.$$
This ends the proof.

\subsection{Proof of Corollary \ref{prop:estimation:R}}\label{sec_cor_R_prob}

Let us introduce the notations
$$
I_n=\sup_{(x,y)\in [1,\infty)\times[0,+\infty)}\left|\mathcal{R}(x,y)-\widehat{\mathcal{R}}_n(x,y)\right|\qquad\text{and}\qquad C=\frac{1}{\lambda_\ast}\left(4e^{-1}\frac{\lambda^\ast }{\lambda_\ast}+1\right).
$$
For any $\e>0$, according to Proposition \ref{prop:maj:R}, we have
\begin{align*}
\prob(I_n\geq\e)&\leq\prob\left(\left\|G-\widehat{G}_n\right\|_\infty + C\left\|\widehat{G}_n\right\|_\infty \left|\lambda-\widehat{\lambda}_n\right|\geq \e\right)\\
&\leq \prob\left(\left\|G-\widehat{G}_n\right\|_\infty\geq\frac{\e}{2}\right)+\prob\left(C\left\|\widehat{G}_n\right\|_\infty \left|\lambda-\widehat{\lambda}_n\right|\geq\frac{\e}{2}\right).
\end{align*}
Let $\eta$ be a positive real. Using the elementary inequality satisfied for any reals $a$ and $b$,
$$
ab\leq \frac{1}{4\eta}a^2+\eta b^2,
$$
we have
\begin{align*}
\prob\left(\|\widehat{G}_n\|_\infty \left|\lambda-\widehat{\lambda}_n\right|\geq\frac{\e}{2C}\right)&\leq \prob\left(\eta\|\widehat{G}_n\|^2_\infty+\frac{1}{4\eta} \left|\lambda-\widehat{\lambda}_n\right|^2\geq\frac{\e}{2C}\right)\\
&\leq \prob\left(\eta\|\widehat{G}_n\|^2_\infty\geq\frac{\e}{4C}\right)+\prob\left(\frac{1}{4\eta} \left|\lambda-\widehat{\lambda}_n\right|^2\geq\frac{\e}{4C}\right).
\end{align*}
Notice that
\begin{align*}
\prob\left(\|\widehat{G}_n\|^2_\infty\geq\frac{\e}{4C\eta}\right)&=\prob\left(\|\widehat{G}_n\|_\infty\geq\sqrt{\frac{\e}{4C\eta}}\right)\\
&\leq \prob\left(\|G\|_\infty+\|\widehat{G}_n-G\|_\infty\geq\sqrt{\frac{\e}{4C\eta}}\right)\\
&\leq\prob\left(\|G\|_\infty\geq\frac12\sqrt{\frac{\e}{4C\eta}}\right)+\prob\left(\|\widehat{G}_n-G\|_\infty\geq\frac12\sqrt{\frac{\e}{4C\eta}}\right).
\end{align*}
With $\eta^\ast=\frac{\e}{16C(\|G\|_\infty+1)^2}$, we have
$$
\prob\left(\|G\|_\infty\geq\frac12\sqrt{\frac{\e}{4C\eta}}\right)=0.
$$
Thus, {we obtain}
$$
\prob\left(\|\widehat{G}_n\|^2_\infty\geq\frac{\e}{4C\eta^\ast}\right)\leq \prob\left(\|\widehat{G}_n-G\|_\infty\geq\frac12\sqrt{\frac{\e}{4C\eta^\ast}}\right).
$$
To sum up,
\begin{equation}\label{eq:rate:convergence}
\prob(I_n\geq\e)\leq \prob\left(\|\widehat{G}_n-G\|_\infty\geq\frac\e2\right)+\prob\left(\|\widehat{G}_n-G\|_\infty\geq\|G\|_\infty+1\right)+\prob\left( \left|\lambda-\widehat{\lambda}_n\right|\geq\frac{\e}{4C(1+\|G\|_\infty)}\right).
\end{equation}
The result follows.

\subsection{Proof of Theorem \ref{prop:estim:r}}

Let us define the operator
\begin{equation}\label{eq:operatorK}
K\,:\,h\mapsto \int_1^{+\infty} h(y)\mathcal{R}(x,y)dy
\end{equation}
on $L^1(1,\infty)$. Let us show that the norm of the operator $K$ on $L^1(1,\infty)$ is less than $1$ if the condition $\frac{\lambda}{\lambda+r}\int_0^1 \frac{G(u)}{u}\,du<1$ is satisfied. Indeed, using Jensen and Fubini's theorems, for any $h\in L^1(1,\infty)$ we have,
$$
\|Kh\|=\int_1^\infty\left|\int_1^\infty h(y)\mathcal{R}(x,y)\,dy\right|\,dx\leq \int_1^\infty\int_1^\infty \left|h(y)\right|\mathcal{R}(x,y)\,dy\,dx\leq \sup_{y\in[1,\infty[}\int_1^\infty \mathcal{R}(x,y)dx\int_1^\infty \left|h(y)\right|dy.
$$
According to Lemma \ref{lem-cond-K}, the above inequalities yield
$$
\|Kh\|\leq \frac{\lambda}{\lambda+r}\int_0^1 \frac{G(u)}{u}\,du\|h\|.
$$
Therefore, under the condition, $\frac{\lambda}{\lambda+r}\int_0^1 \frac{G(u)}{u}\,du<1$, we get $\|K\|<1$. One may then rewrite equation (\ref{eq:int:r}) as the Fredholm equation
$$
p-Kp=s,
$$
where $s(x)=\int_0^1\mathcal{R}(x,y)\,dy$. This equation has obviously a unique solution since $\|K\|<1$. Notice that $\|p\|<\infty$ since one may write $p=\sum_{k=0}^\infty K^ks$. The following proposition precises the relations between $\widehat{K}_n$ with $K$ and $\widehat{s}_n$ with $s$.
\begin{proposition}\label{K-hatK}
The estimation $\widehat{K}_n$ and $\widehat{s}_n$ converge toward $K$ and $s$ respectively in probability. For any $\e>0$,
$$
\lim_{n\to\infty}\prob\left(\|\widehat{K}_n-K\|\geq \e\right)=0,\,\lim_{n\to\infty}\prob\left(\|\widehat{s}_n-s\|\geq \e\right)=0.
$$
\end{proposition}
\proof
First, let us notice that for any $h\in L^1(1,\infty)$, we have
\begin{align*}
\|(K-\widehat{K}_n)h\|&=\int_1^\infty\left|\int_1^\infty h(y)(\mathcal{R}(x,y)-\widehat{\mathcal{R}}_n(x,y))\,dy\right|\,dx\\
&\leq \int_1^\infty\int_1^\infty\left| h(y)\right|\left|\mathcal{R}(x,y)-\widehat{\mathcal{R}}_n(x,y)\right|\,dy\,dx\\
&\leq \|h\|\sup_{y\geq1}\int_1^\infty\left|\mathcal{R}(x,y)-\widehat{\mathcal{R}}_n(x,y)\right|\,dx.
\end{align*}
Therefore, 
$$
\|K-\widehat{K}_n\|\leq \sup_{y\geq1}\int_1^\infty\left|\mathcal{R}(x,y)-\widehat{\mathcal{R}}_n(x,y)\right|\,dx
$$
$\prob$-a.s. Then, using Lemma \ref{hatR-specialnorm}, almost-surely we have
$$
\|K-\widehat{K}_n\|\leq\frac{\lambda}{\lambda+r}\int_0^1 |G(u)-\widehat{G}_n(u)|\,u^{-1}\,du+\lambda^\ast\left(4e^{-1}\frac{1}{\lambda^2_\ast }+\frac{1}{\lambda_\ast+r}\right)\int_0^1 |\widehat{G}_n(u)|\,u^{-1}\,du\,|\lambda-\widehat{\lambda}_n|.
$$
Then, using Assumption \ref{hyp:G}, $(C^\lambda_{1,2})$ and $(C^G_2)$, the convergence in probability of $\widehat{K}_n$ towards $K$ follows. The proof of the convergence in probability of $\widehat{s}_n$ towards $s$ in probability is quite similar.
\hfill$\Box$

\noindent We deduce easily from the above proposition that for any $\e>0$,
$$
\lim_{n\to\infty}\prob\left(\|(\widehat{K}_n-K)r\|\geq \e\right)=0.
$$
Let us choose $\eta>0$ and $\e>0$ such that $\e<1-\|K\|$. We define
$$
\Omega_{n}=\left\{\omega\in\Omega ~ ; ~ \|\widehat{K}_n(\omega)\|<1-\e, \|(K-\widehat{K}_n(\omega))r\|\leq\frac{\e^2}{4}, \|s-\widehat{s}_{n}(\omega)\|\leq\frac{\e^2}{4}\right\}.
$$
According to Proposition \ref{K-hatK}, there exists $N$ such that for all $n\geq N$,
\begin{align*}
\prob(\Omega\setminus\Omega_{n})&\,\leq\,\prob\left(\|\widehat{K}_n\|\geq1-\e\right)+\prob\left(\|(K-\widehat{K}_n)r\|>\frac{\e^2}{4}\right)+\prob\left(\|s-\widehat{s}_{n}\|>\frac{\e^2}{4}\right)\\
&\,\leq\,\frac{\eta}{3}+\frac{\eta}{3}+\frac{\eta}{3}\,=\,\eta
\end{align*}
From (\ref{eq:int:r:hat}), $\widehat{p}_n$ satisfies almost-surely the equation,
$$
\widehat{p}_n=\widehat{s}_n+\widehat{K}_n\widehat{p}_n.
$$
Therefore, we also have $\widehat{p}_n=\sum_{k=0}^\infty\widehat{K}^k_n\widehat{s}_{n}$.
We split the difference $p-\widehat{p}_{n,m}$ using the quantity $\widehat{p}_n$,
$$p-\widehat{p}_{n,m}\,=\,p-\widehat{p}_{n}~+~\widehat{p}_{n}-\widehat{p}_{n,m} .$$ We begin to bound $p-\widehat{p}_{n}$ on $\Omega_{n}$. For $n\geq N$, on $\Omega_{n}$,
\begin{align*}
\|p-\widehat{p}_n\|&\leq \|(s-\widehat{s}_{n})p\|+\|(K-\widehat{K}_n)p\|+\|\widehat{K}_n\|\|p-\widehat{p}_n\|\\
&\leq \frac{\e^2}{4}+\frac{\e^2}{4}+(1-\e)\|p-\widehat{p}_n\|.
\end{align*}
An elementary re-arranging yields $\|p-\widehat{p}_n\|\leq \frac{\e}{2}$ on $\Omega_{n}$. It remains to consider the difference $\widehat{p}_n-\widehat{p}_{n,m}$. By definition, we have
$$
\widehat{p}_n-\widehat{p}_{n,m}=\sum_{k=m+1}^\infty\widehat{K}^k_n\widehat{s}_{n}.
$$
Therefore, for $n\geq N$, on $\Omega_{n}$, {we have}
$$
\|\widehat{p}_n-\widehat{p}_{n,m}\|~\leq~\|\widehat{s}_{n}\|\frac{\|\widehat{K}_n\|^{m+1}}{1-\|\widehat{K}_n\|}~\leq~\left(\frac{\e^2}{2}+\|s\|\right)\frac{(1-\e)^{m+1}}{\e}~\leq~\frac{\e}{2}
$$
for $m\geq M$ with $M$ {large} enough. Therefore, for $n\geq N$ and $m\geq M$, on  $\Omega_{n}$,
$$
\|p-\widehat{p}_{n,m}\|~\leq~\|p-\widehat{p}_{n}\|+\|\widehat{p}_n-\widehat{p}_{n,m}\|~\leq~\frac{\e}{2}+\frac{\e}{2}~=~\e.
$$
This concludes the proof.

\subsection{Proof of Theorem \ref{prop:estim:t}}

All the ingredients for this proof are in fact already present in the proof of Theorem \ref{prop:estim:r}. Nevertheless, let us give some details. As in the previous section, let us choose $\eta>0$ and $\e>0$ such that $\e<1-\|K\|$. We define
$$
\Omega_{n}=\left\{\omega\in\Omega~;~\|\widehat{K}_n(\omega)\|<1-\e, \|(K-\widehat{K}_n(\omega))\|\leq\frac{\e^2}{2\|s\|}, \|s-\widehat{s}_{n}(\omega)\|\leq\frac{\e}{2}\right\}.
$$
According to Proposition \ref{K-hatK}, there exists $N$ such that, for all $n\geq N$,
\begin{align*}
\prob(\Omega\setminus\Omega_{n})&\leq\eta.
\end{align*}
For $m=1$ we have,
$$
\|\widehat{t}_{1,n}-t_1\|=\|\widehat{s}_{n}-s\|\leq\frac{\e}{2}
$$
on $\Omega_n$ and the result follows. Now for $m\geq2$, one may write
$$
\widehat{t}_{m,n}-t_m=\widehat{K}_n(\widehat{t}_{m-1,n}-t_{m-1})+(\widehat{K}_n-K)t_{m-1}.
$$
Notice that on $\Omega_n$, for $m\geq1$,
$$
\|t_m\|=\|K^{m-1}s\|\leq\|K\|^{m-1}\|s\|\leq (1-\e)^{m-1}\|s\|.
$$
Then, for $m\geq2$, {we write}
\begin{align*}
\|\widehat{t}_{m,n}-t_m\|&\leq\|\widehat{K}_n(\widehat{t}_{m-1,n}-t_{m-1})\|+\|(\widehat{K}_n-K)t_{m-1}\|\\
&\leq\|\widehat{K}_n\|\|\widehat{t}_{m-1,n}-t_{m-1}\|+\|(\widehat{K}_n-K)\|\|t_{m-1}\|\\
&\leq(1-\e)\|\widehat{t}_{m-1,n}-t_{m-1}\|+\frac{\e^2}{2\|s\|}(1-\e)^{m-2}\|s\|.
\end{align*}
A straightforward recursion gives, always for $m\geq2$ and on $\Omega_n$,
$$
\|\widehat{t}_{m,n}-t_m\|\leq(1-\e)^{m-1}\|\widehat{t}_{1,n}-t_1\|+\frac{\e^2}{2}(1-\e)^{m-2}\sum_{k=0}^{m-2}(1-\e)^k.
$$
Therefore, for any $m\geq2$, on $\Omega_n$, {we obtain}
$$
\|\widehat{t}_{m,n}-t_m\|~\leq~(1-\e)^{m-1}\frac{\e}{2}+\frac{\e^2}{2}(1-\e)^{m-2}\frac{1-(1-\e)^{m-1}}{\e}~\leq~\e.
$$
\section{Discussion on the condition $(C_2^G)$}
\label{sec:discu}

Here, we propose to show that the Parzen-Rosenblatt estimator $\widehat{G}_n^{PR}$ of the density $G$, defined by
$$\forall\,x\in[0,1],~\widehat{G}_n^{PR}(x) = \frac{1}{n h_n} \sum_{i=1}^n \mathbb{K}\left(\frac{Y_i-x}{h_n}\right),$$
where $\mathbb{K}$ is a kernel function and the bandwidth sequence $(h_n)_{n\geq1}$ tends to $0$ as $n$ goes to infinity, satisfies the condition $(C_2^G)$ under the following assumption on the density of interest.

\begin{hypothese}\label{hyp:G2}
We assume that there exists a real number $\epsilon_1>0$ such that, for any $0\leq x<\epsilon_1$, $G(x)=0$. In addition, we suppose that $G$ is in the H\"older class $\Sigma(\beta,L)$ \citep[Definition 1.2]{TSY}.
\end{hypothese}

\begin{remarque}
Notice that Assumption \ref{hyp:G2} does not hold true for the example in Section \ref{sec:simulations} where $G(u)=11u^{10}$ for $u\in[0,1]$. However, the numerical illustrations show that our results still apply in this case. It means that Assumption \ref{hyp:G2} is certainly non optimal and that some weaker assumptions on $G$, as being close enough to zero near zero, may be sufficient for our results to apply.
\end{remarque}

\noindent
For any $x$, we define the mean squared error of $\widehat{G}_n^{PR}(x)$ by
$$MSE(x)=\E\left[ \left(\widehat{G}_n^{PR}(x)-G(x)\right)^2\right].$$
By \citep[equation $(1.4)$]{TSY}, we have the following bias-variance decomposition
$$MSE(x) = b^2(x)\,+\,V(x),$$
where, with \citep[equation $(1.6)$]{TSY},
$$b(x) = \E\left[\widehat{G}_n^{PR}(x)\right]-G(x)  \qquad\text{and}\qquad V(x)=\frac{1}{nh_n^2}\E\left[\mathbb{K}^2\left(\frac{Y_1-x}{h_n}\right)\right].$$

\noindent
In the sequel, we assume that the chosen kernel function $\mathbb{K}$ has a bounded support. As a consequence, for $n$ large enough and some $\epsilon_2>0$, $\mathbb{K}\left(\frac{y-x}{h_n}\right)=0$ for any $x<\epsilon_2$ and $y\geq\epsilon_1$. Thus,
\begin{equation}\label{eq:pr:001}
\int_0^1\frac{\sqrt{V(x)}}{x}dx \,=\, \frac{1}{\sqrt{n} h_n} \int_{\epsilon_1}^1 \frac{1}{x}\left[\int_{\epsilon_2}^1 G(y)\,\mathbb{K}^2\left(\frac{y-x}{h_n}\right) dx\right]^{1/2} dy \,\leq\, -\frac{\|\mathbb{K}\|_{\infty} \log(\epsilon_1)}{\sqrt{n} h_n} .
\end{equation}
\noindent
In addition, $b(x)=0$ for any $x<\epsilon_1\wedge\epsilon_2$. Therefore, by virtue of \citep[Proposition 1.2]{TSY},
\begin{equation}\label{eq:pr:002}
\int_0^1 \frac{b(x)}{x}dx \,=\,\int_{\epsilon_1\wedge\epsilon_2}^1 \frac{b(x)}{x}dx \,\leq\,C_1\,h_n^\beta ,\end{equation}
for some positive number $C_1$, whenever $\mathbb{K}$ is a kernel of order $l=\lfloor\beta\rfloor$ (see \citep[Definition 1.3]{TSY}) satisfying
$$\int |u|^\beta \mathbb{K}(u) du\,<\,\infty.$$
Finally, by $(\ref{eq:pr:001})$ and $(\ref{eq:pr:002})$, we have
\begin{eqnarray*}
\E\left[\int_0^1 \frac{\big|\widehat{G}_n^{PR}(x)-G(x)\big|}{x} dx \right]&\leq& \int_0^1\frac{\sqrt{MSE(x)}}{x} dx\\
&\leq&\int_0^1\frac{b(x)}{x} dx \,+\,\int_0^1 \frac{\sqrt{V(x)}}{x} dx \\
&\leq& C_2\left( \frac{1}{\sqrt{n}h_n} + h_n^\beta\right),
\end{eqnarray*}
for some constant $C_2$. We conclude that the $L^1$-norm vanishes when $n$ tends to infinity if the bandwidth is such that $\sqrt{n}h_n\to0$. Therefore, the convergence in probability $(C_2^G)$ holds under this condition.

\bibliographystyle{plainnat}
\bibliography{tex_main}

\begin{thebibliography}{21}
\providecommand{\natexlab}[1]{#1}
\providecommand{\url}[1]{\texttt{#1}}
\expandafter\ifx\csname urlstyle\endcsname\relax
  \providecommand{\doi}[1]{doi: #1}\else
  \providecommand{\doi}{doi: \begingroup \urlstyle{rm}\Url}\fi

\bibitem[Aza\"\i{}s(2014)]{AzaisESAIM14}
Romain Aza\"\i{}s.
\newblock A recursive nonparametric estimator for the transition kernel of a
  piecewise-deterministic {M}arkov process.
\newblock \emph{To appear in ESAIM: Probability and Statistics}, 2014.

\bibitem[Aza\"\i{}s et~al.(2014)Aza\"\i{}s, Dufour, and
  G\'egout-Petit]{AzaisSJS14}
Romain Aza\"\i{}s, Fran\c{c}ois Dufour, and Anne G\'egout-Petit.
\newblock Nonparametric estimation of the conditional distribution of the
  inter-jumping times for piecewise-deterministic {M}arkov processes.
\newblock \emph{To appear in Scandinavian Journal of Statistics}, 2014.

\bibitem[Bena\"\i{}m et~al.(2014)Bena\"\i{}m, Le~Borgne, Malrieu, and
  Zitt]{benaim2014}
Michel Bena\"\i{}m, St\'ephane Le~Borgne, Florent Malrieu, and Pierre-Andr\'e
  Zitt.
\newblock On the stability of planar randomly switched systems.
\newblock \emph{The Annals of Applied Probability}, 24\penalty0 (1):\penalty0
  292--311, 02 2014.

\bibitem[Brandejsky et~al.(2013)Brandejsky, De~Saporta, and
  Dufour]{brandejsky:spa}
Adrien Brandejsky, Beno{\^\i}te De~Saporta, and Fran{\c c}ois Dufour.
\newblock {Optimal stopping for partially observed piecewise-deterministic
  Markov processes}.
\newblock \emph{Stochastic Processes and their Applications}, 123:\penalty0
  3201--3238, 2013.

\bibitem[Buckwar and Riedler(2011)]{buckwar2011exact}
Evelyn Buckwar and Martin~G. Riedler.
\newblock An exact stochastic hybrid model of excitable membranes including
  spatio-temporal evolution.
\newblock \emph{Journal of mathematical biology}, 63\penalty0 (6):\penalty0
  1051--1093, 2011.

\bibitem[Chiquet and Limnios(2008)]{ChiquetLimnios}
Julien Chiquet and Nikolaos Limnios.
\newblock A method to compute the transition function of a piecewise
  deterministic {M}arkov process with application to reliability.
\newblock \emph{Statist. Probab. Lett.}, 78\penalty0 (12):\penalty0 1397--1403,
  2008.

\bibitem[Cloez and Hairer(2014)]{cloezhairer}
Bertrand Cloez and Martin Hairer.
\newblock Exponential ergodicity for markov processes with random switching.
\newblock \emph{To appear in Bernoulli}, 2014.

\bibitem[Costa and Dufour(2008)]{costa:siam}
Oswaldo Costa and Fran{\c c}ois Dufour.
\newblock {Stability and Ergodicity of Piecewise Deterministic Markov
  Processes}.
\newblock \emph{SIAM Journal on Control and Optimization}, 47\penalty0
  (2):\penalty0 1053--1077, 2008.

\bibitem[Costa and Dufour(2013)]{CD:book}
Oswaldo Costa and Fran{\c{c}}ois Dufour.
\newblock \emph{Continuous Average Control of Piecewise Deterministic Markov
  Processes}.
\newblock SpringerBriefs in Mathematics. Springer Verlag, 2013.

\bibitem[Davis(1993)]{DAVIS}
Mark H.~A. Davis.
\newblock \emph{Markov models and optimization}, volume~49 of \emph{Monographs
  on Statistics and Applied Probability}.
\newblock Chapman \& Hall, London, 1993.

\bibitem[De~Saporta et~al.(2012)De~Saporta, Dufour, Zhang, and Elegbede]{DeS}
Beno{\^\i}te De~Saporta, Fran\c{c}ois Dufour, Huilong Zhang, and Charles
  Elegbede.
\newblock Optimal stopping for the predictive maintenance of a structure
  subject to corrosion.
\newblock \emph{Journal of Risk and Reliability}, 226 (2):\penalty0 169--181,
  2012.

\bibitem[Doumic et~al.(2011)Doumic, Hoffmann, Reynaud-Bouret, and
  Rivoirard]{Doumic12}
Marie Doumic, Marc Hoffmann, Patricia Reynaud-Bouret, and Vincent Rivoirard.
\newblock Nonparametric estimation of the division rate of a size-structured
  population.
\newblock \emph{SIAM Journal on Numerical Analysis}, 50\penalty0 (2):\penalty0
  925--950, 2011.

\bibitem[Doumic et~al.(2014)Doumic, Hoffmann, Krell, and Robert]{Doumic11}
Marie Doumic, Marc Hoffmann, Nathalie Krell, and Lydia Robert.
\newblock Statistical estimation of a growth-fragmentation model observed on a
  genealogical tree.
\newblock \emph{To appear in Bernoulli}, 2014.

\bibitem[Genadot and Thieullen(2012)]{genadot2012averaging}
Alexandre Genadot and Mich{\`e}le Thieullen.
\newblock Averaging for a fully coupled piecewise-deterministic markov process
  in infinite dimensions.
\newblock \emph{Advances in Applied Probability}, 44\penalty0 (3):\penalty0
  749--773, 2012.

\bibitem[Jacobsen(2006)]{JACOB}
Martin Jacobsen.
\newblock \emph{Point process theory and applications : marked point and
  piecewise deterministic processes}.
\newblock Probability and its applications. Birkh\"{a}user, Boston (Mass.),
  Basel, Berlin, 2006.

\bibitem[Kovacevic and Pflug(2011)]{TJRI}
Raimund~M. Kovacevic and Georg~Ch. Pflug.
\newblock Does insurance help to escape the poverty trap? -- {A} ruin theoretic
  approach.
\newblock \emph{The Journal of Risk and Insurance}, 78\penalty0 (4):\penalty0
  1003--1027, 2011.

\bibitem[Murray(2002)]{MURRAY}
James~D. Murray.
\newblock \emph{Mathematical Biology I: An Introduction}, volume~17 of
  \emph{Interdisciplinary Applied Mathematics}.
\newblock Springer, New York, 2002.

\bibitem[Niethammer and Schempp(1970)]{NIETHAMMER}
Wilhelm Niethammer and Walter Schempp.
\newblock On the construction of iteration methods for linear equations in
  banach spaces by summation methods.
\newblock \emph{aequationes mathematicae}, 5\penalty0 (1):\penalty0 124--125,
  1970.

\bibitem[Riedler et~al.(2012)Riedler, Thieullen, and Wainrib]{riedler2012limit}
Martin~G. Riedler, Mich{\`e}le Thieullen, and Gilles Wainrib.
\newblock Limit theorems for infinite-dimensional piecewise deterministic
  markov processes. applications to stochastic excitable membrane models.
\newblock \emph{Electron. J. probab}, 17\penalty0 (55):\penalty0 1--48, 2012.

\bibitem[Tsybakov(2008)]{TSY}
Alexandre~B. Tsybakov.
\newblock \emph{Introduction to Nonparametric Estimation}.
\newblock Springer Publishing Company, Incorporated, 1st edition, 2008.

\bibitem[Wied and Wei\ss{}bach(2010)]{Wied2010Consistency}
Dominik Wied and Rafael Wei\ss{}bach.
\newblock Consistency of the kernel density estimator - a survey.
\newblock \emph{Statistical Papers}, 53\penalty0 (1):\penalty0 1--21, 2010.

\end{thebibliography}

\begin{figure}[p]
	\centering
	\includegraphics[height=5cm,width=4.1cm]{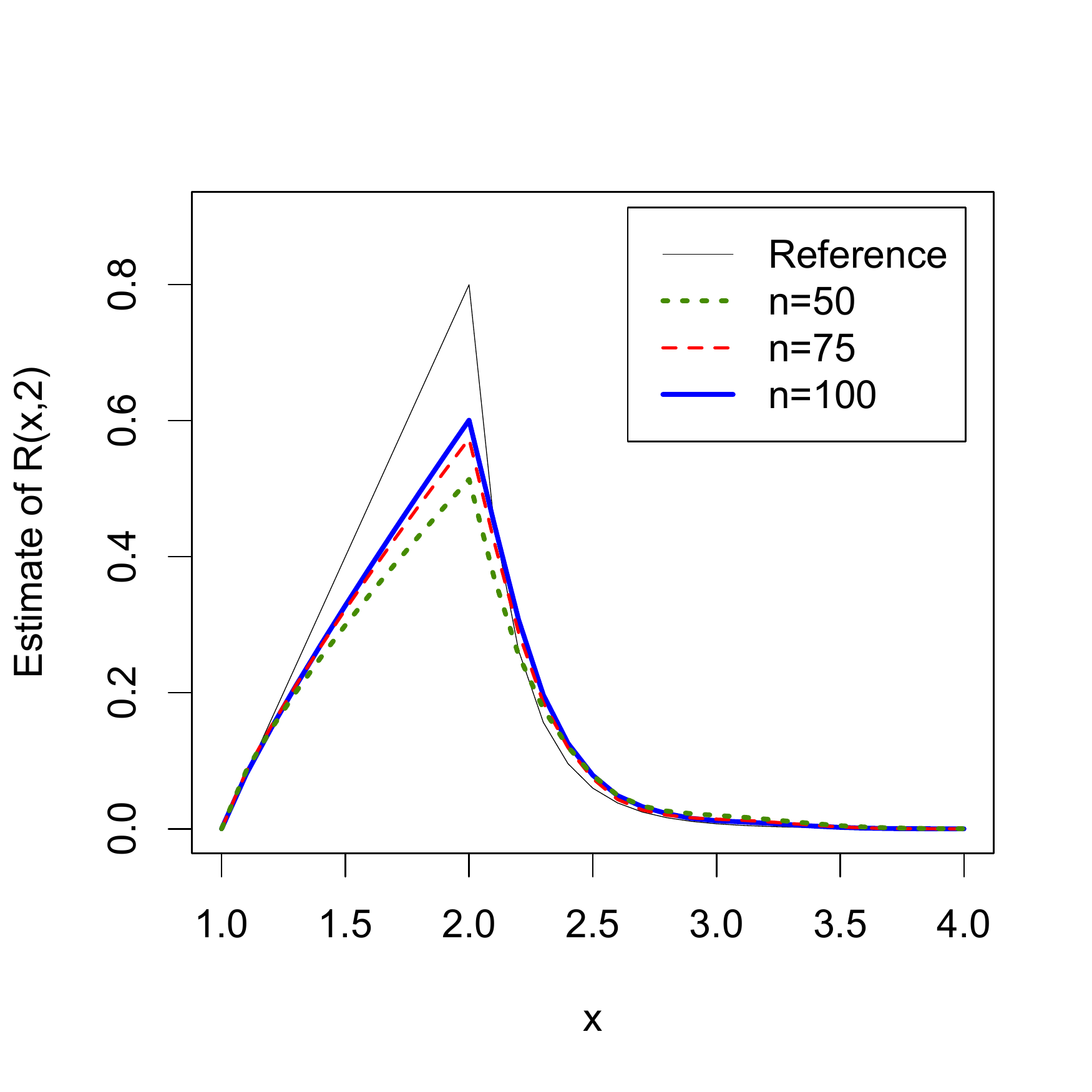}
	\includegraphics[height=5cm,width=4.1cm]{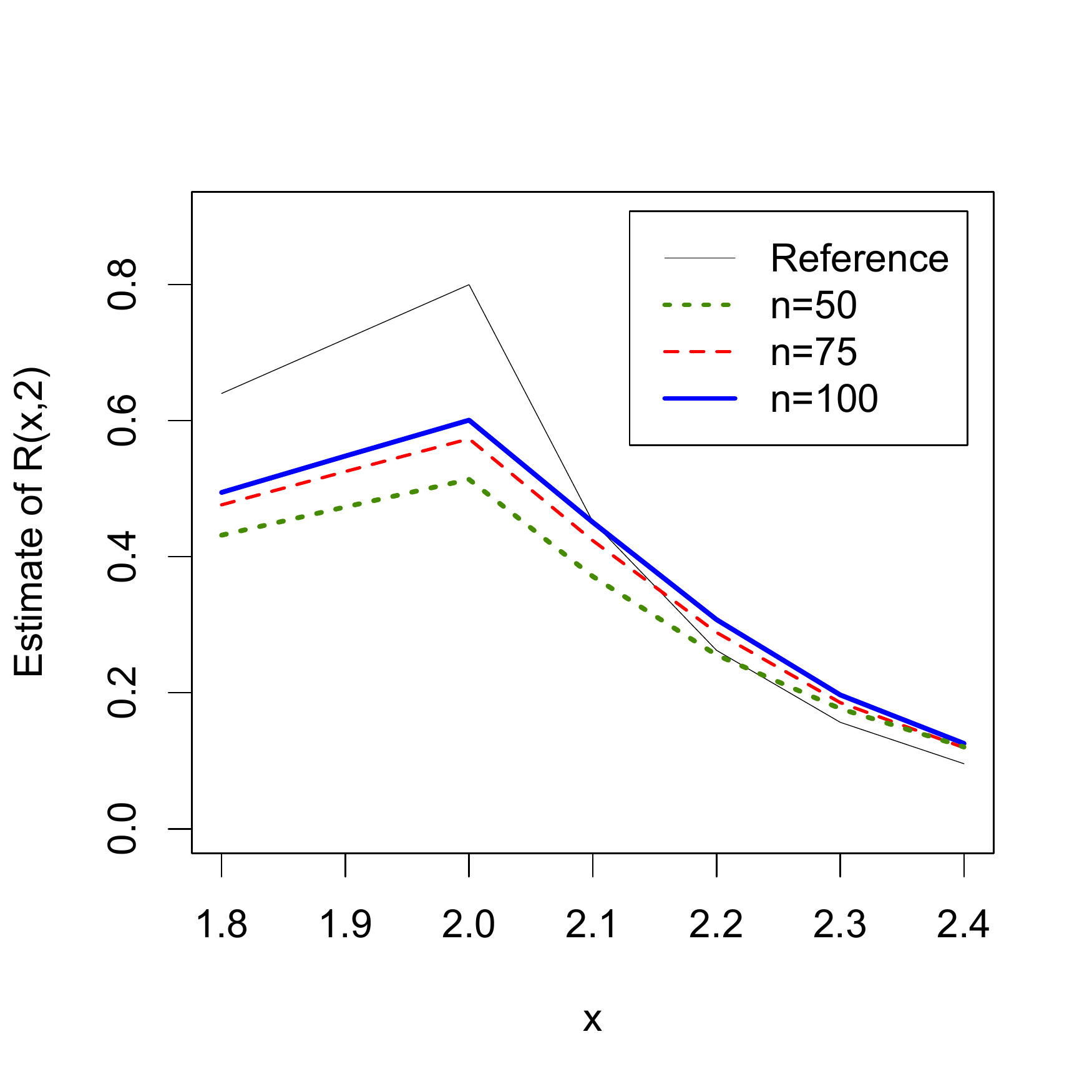}
	\includegraphics[height=4.95cm,width=8cm]{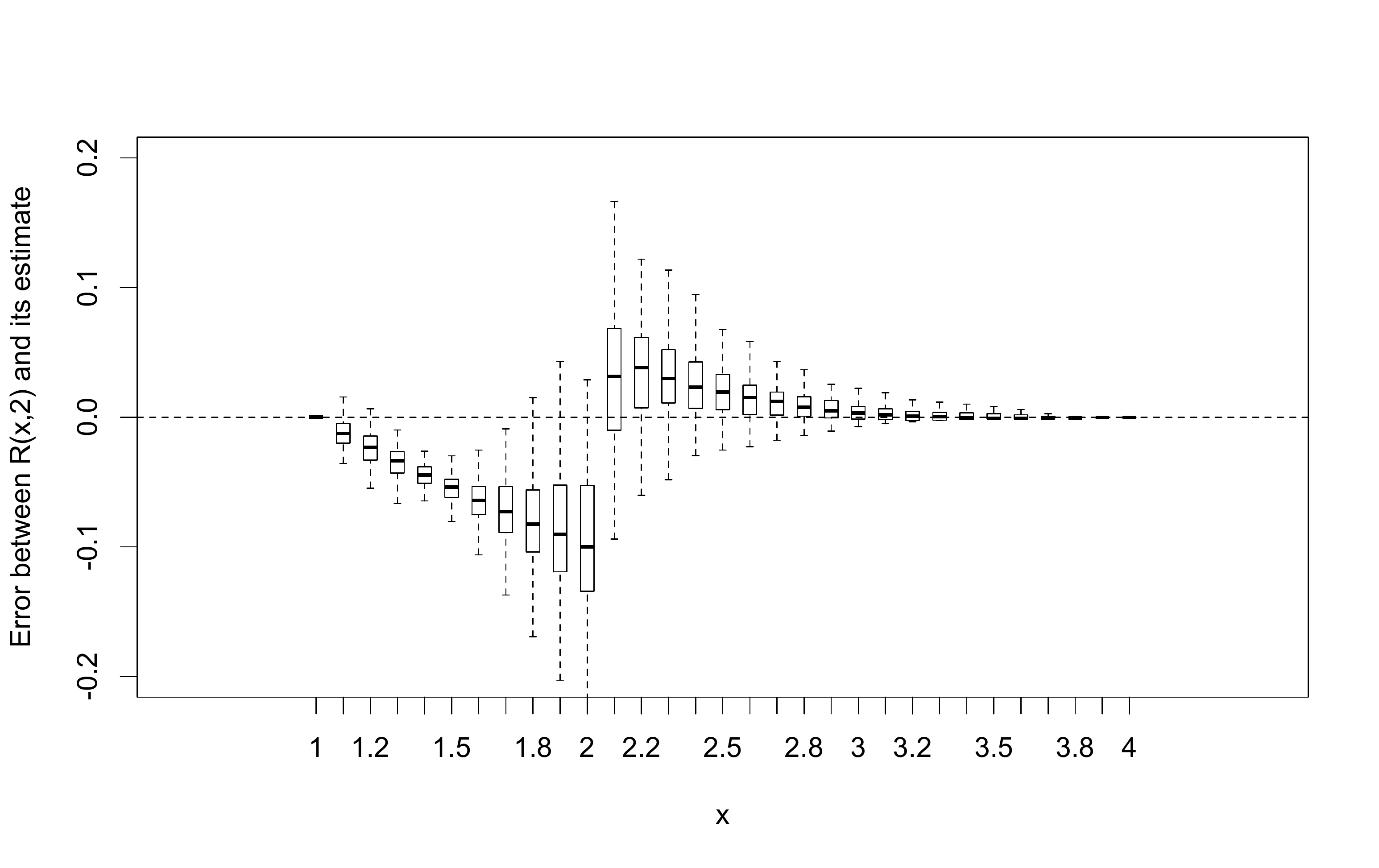}	
\caption{{The figure displays the reference curve $\mathcal{R}(\cdot,2)$ and its estimates from the observation of $n=50,\,75$ or $100$ random loss events (left) with a zoom around $\mathcal{R}(2,2)$ (center), and the pointwise error on $100$ replicates between $\mathcal{R}(\cdot,2)$ and $\widehat{\mathcal{R}}_{100}(\cdot,2)$ (right)}.}
\label{fig:Rx2}
\end{figure}

\begin{figure}[p]
	\centering
	\includegraphics[height=5cm,width=4.1cm]{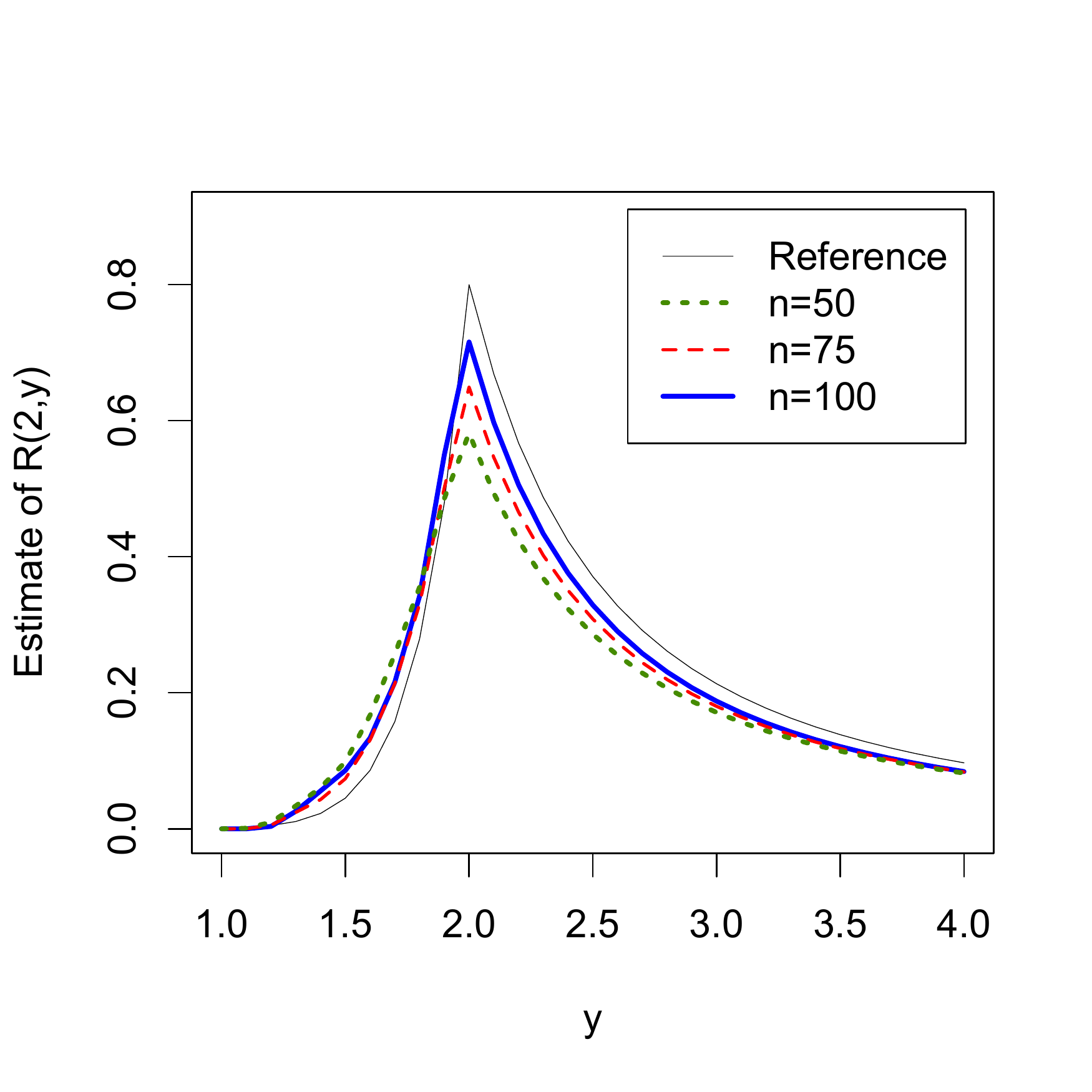}
	\includegraphics[height=5cm,width=4.1cm]{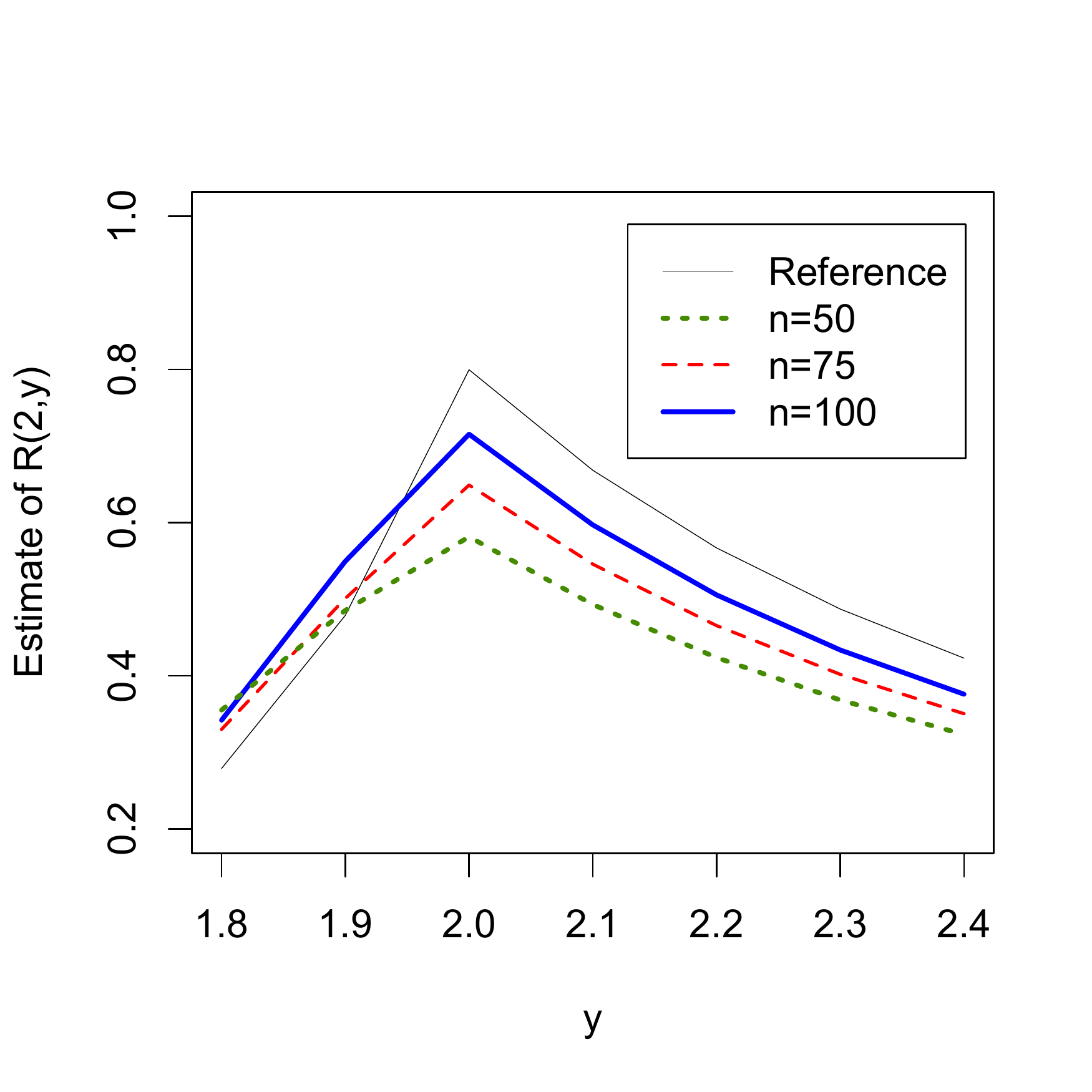}
	\includegraphics[height=4.95cm,width=8cm]{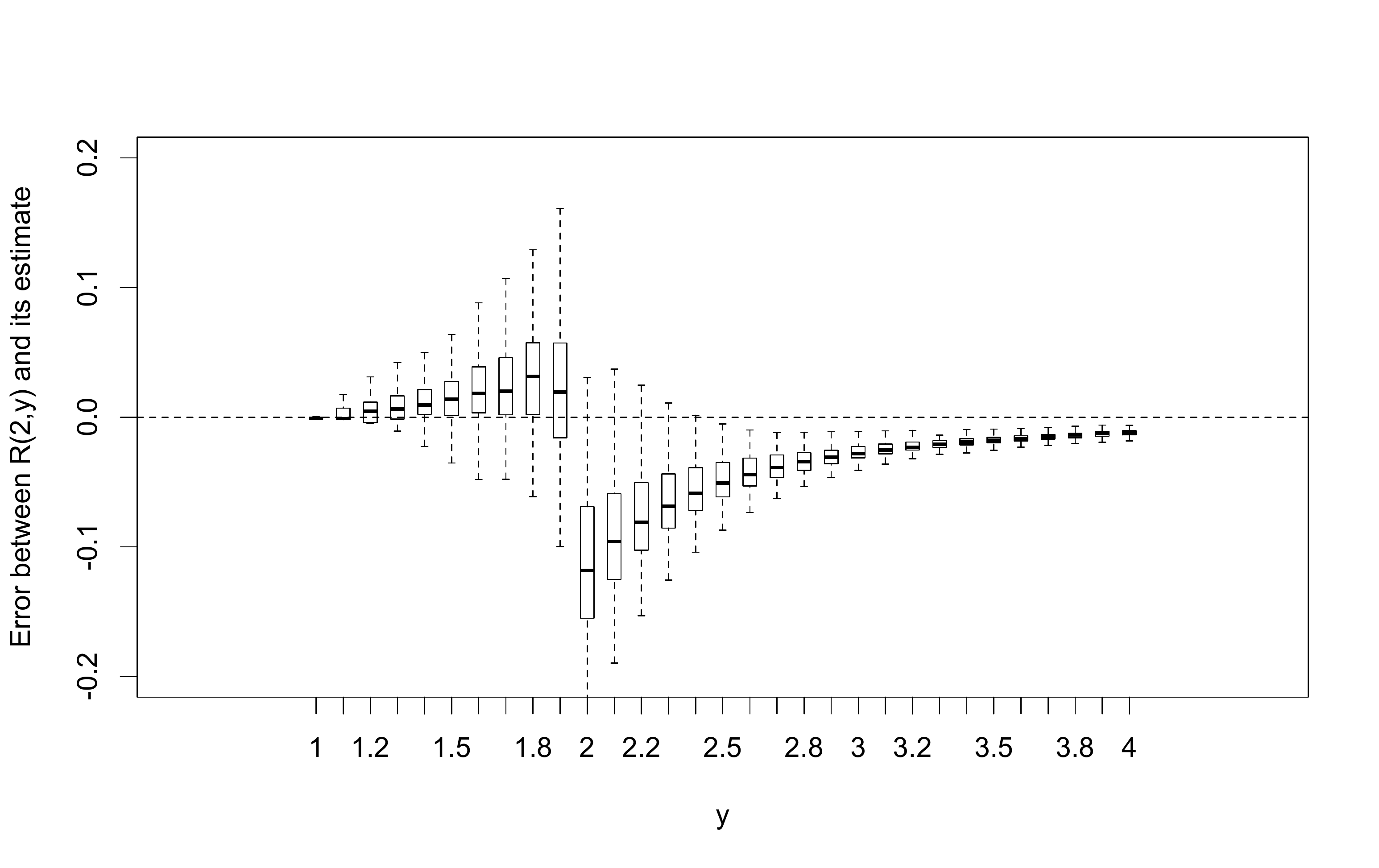}	
\caption{{The figure displays the reference curve $\mathcal{R}(2,\cdot)$ and its estimates from the observation of $n=50,\,75$ or $100$ random loss events (left) with a zoom around $\mathcal{R}(2,2)$ (center), and the pointwise error on $100$ replicates between $\mathcal{R}(2,\cdot)$ and $\widehat{\mathcal{R}}_{100}(2,\cdot)$ (right)}.}
	\label{fig:R2y}
\end{figure}

\begin{figure}[p]
	\centering
	\includegraphics[height=5cm]{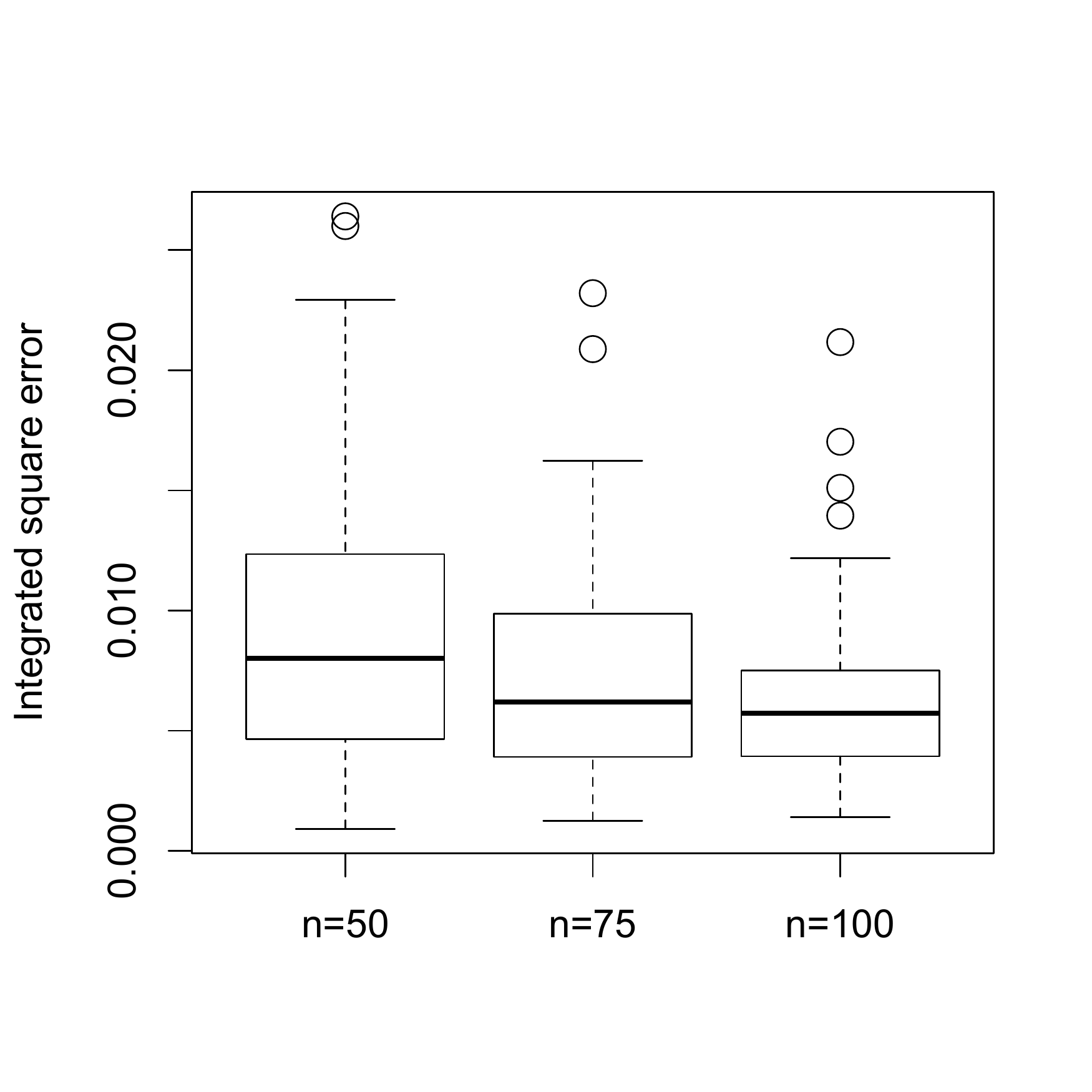}
	\qquad\qquad
	\includegraphics[height=5cm]{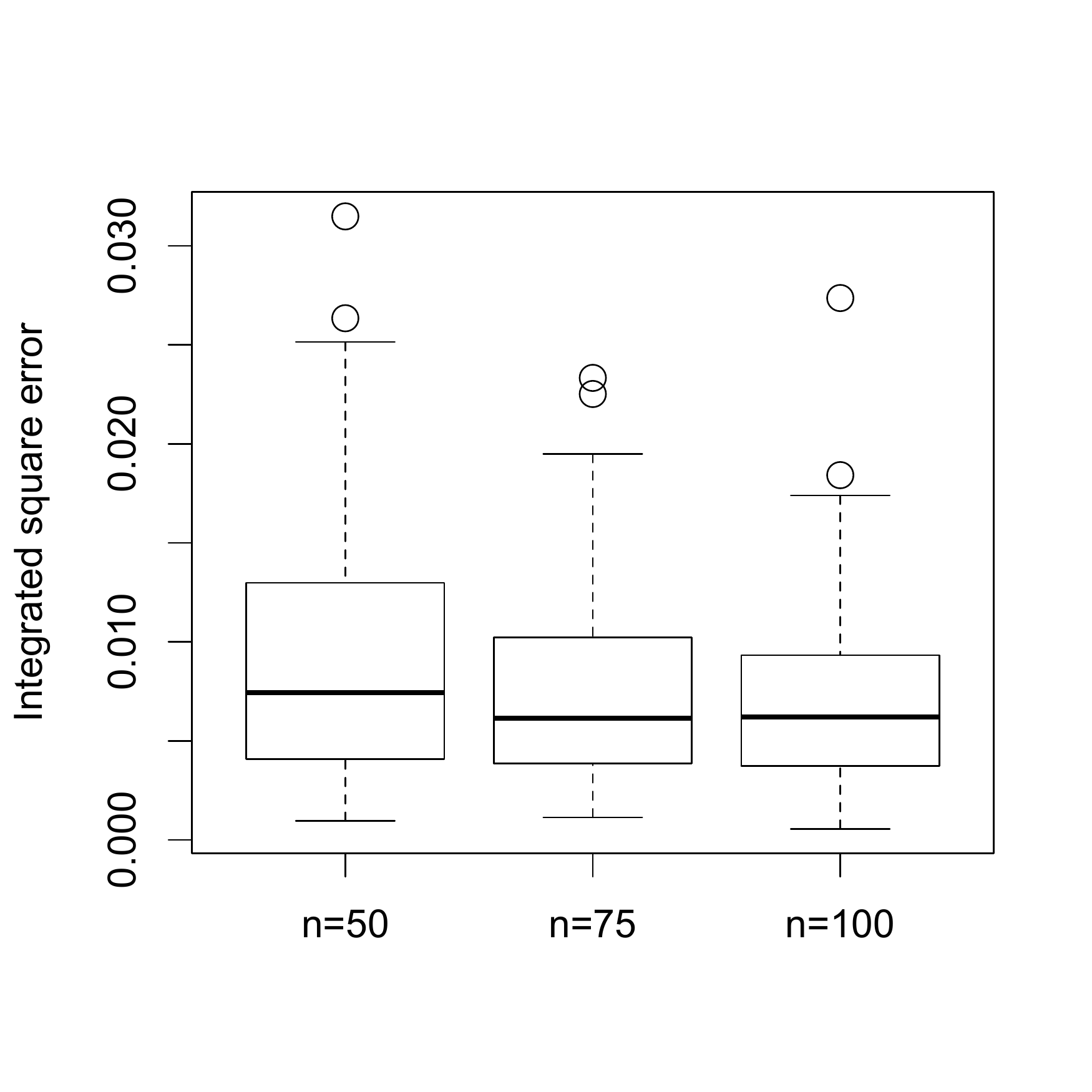}
\caption{{The figure displays the integrated square error on $100$ replicates between $\mathcal{R}(\cdot,2)$ and its estimate (left), and between $\mathcal{R}(2,\cdot)$ and its estimate (right), from the observation of $n=50,\,75$ or $100$ random loss events.}}
\label{fig:R:ISE}
\end{figure}

\begin{figure}[p]
	\centering
	\includegraphics[height=6.2cm]{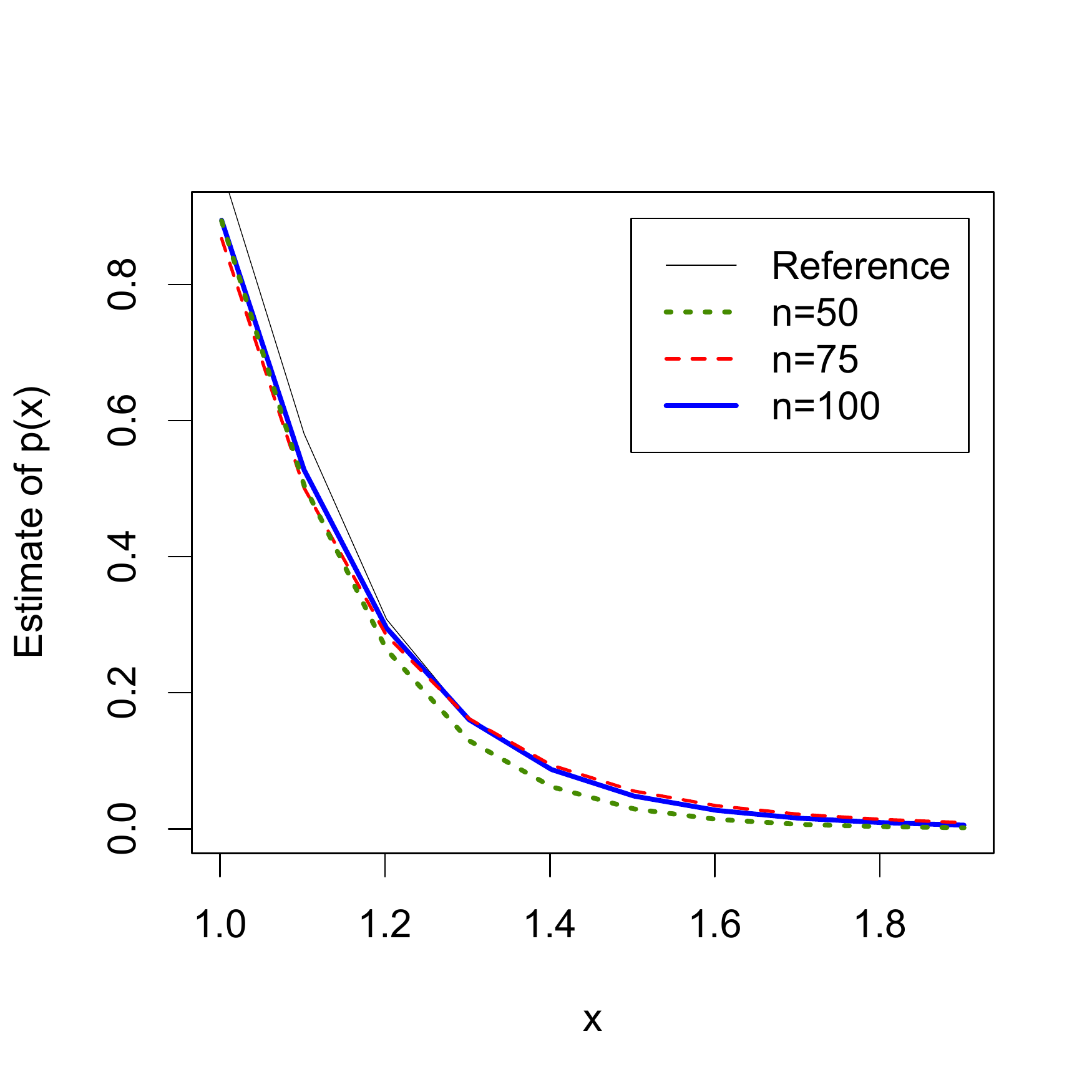}
	\includegraphics[height=6cm]{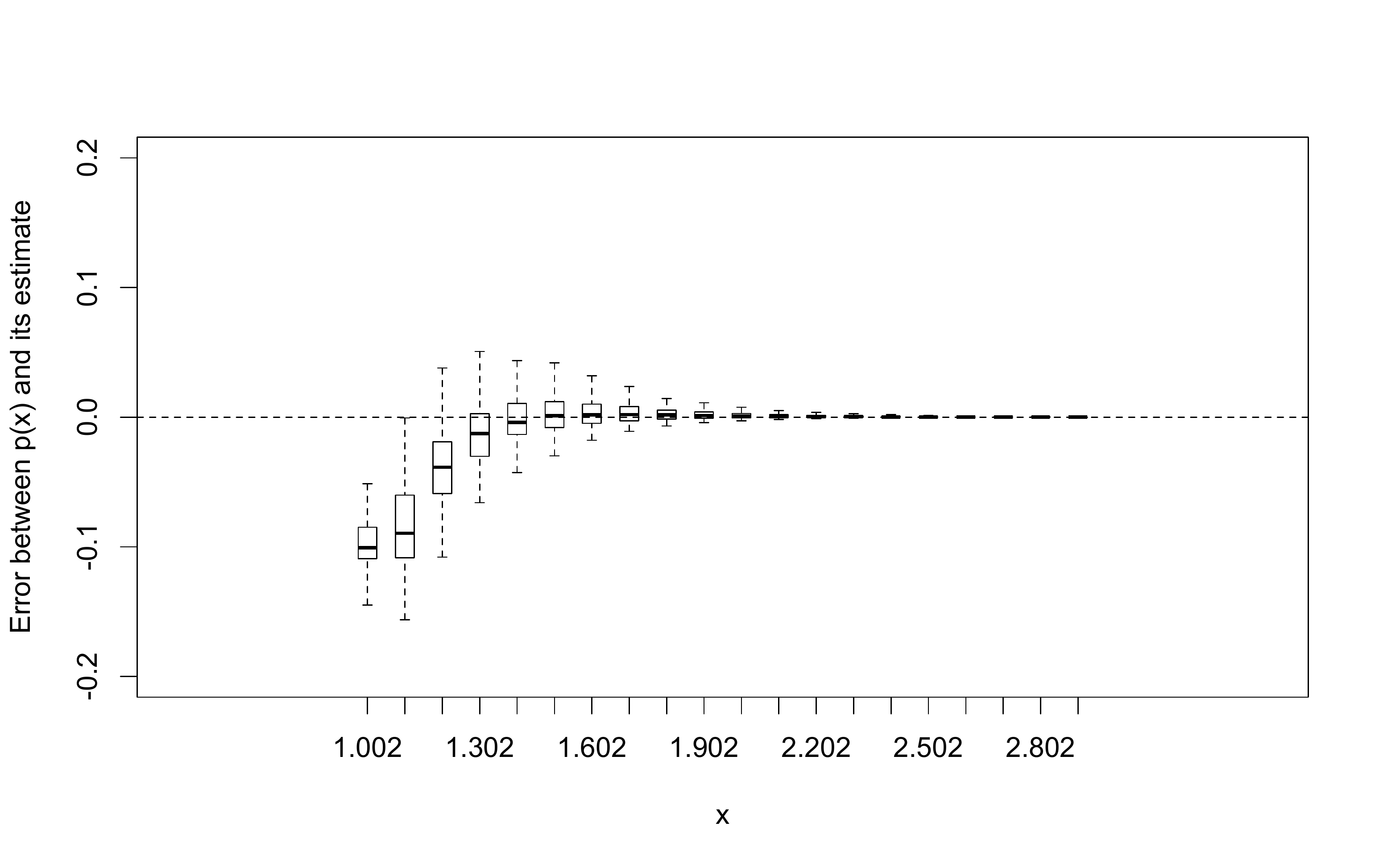}
	\caption{{The figure displays} the absorption probability $p$ (approximated by $p_m$) and its estimates $\widehat{p}_{n,m}$ from the observation of $n=50,\,75$ or $100$ random loss events and for $m=10$ iterations of the estimated kernel $\widehat{K}_n$ (left), and {the associated pointwise error on $100$ replicates from $n=100$ random loss events (right)}.}
	\label{fig:p:curve}
\end{figure}

\begin{figure}[p]
	\centering
	\includegraphics[height=6cm]{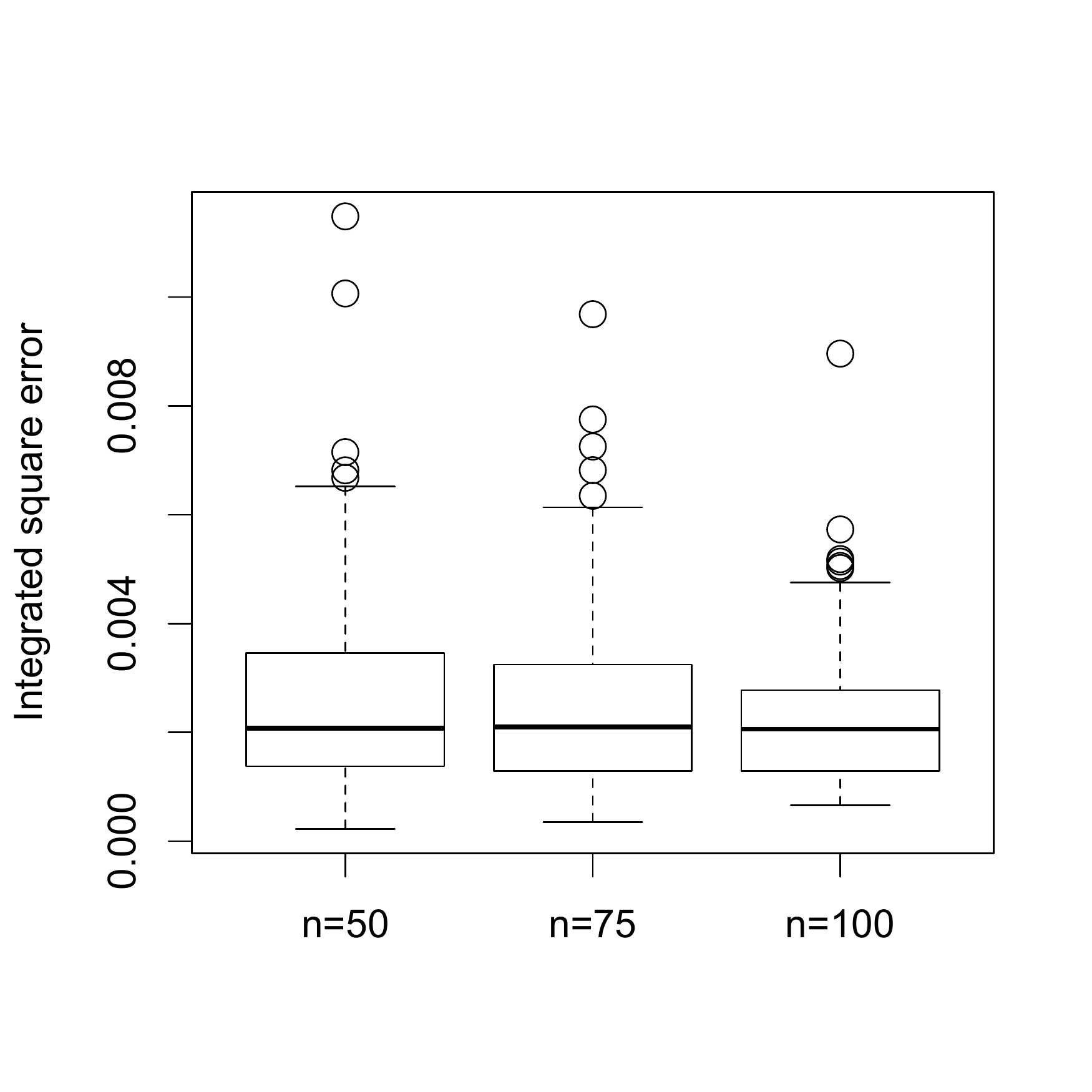}
	\caption{{The figure displays} the integrated square error on $100$ replicates between the absorption probability $p$ (approximated by $p_m$) and its estimates $\widehat{p}_{n,m}$ from the observation of $n=50,\,75$ or $100$ random loss events and for $m=10$ iterations of the estimated kernel $\widehat{K}_n$.}
	\label{fig:p:ise}
\end{figure}

\begin{figure}[p]
	\centering
	\includegraphics[height=5cm]{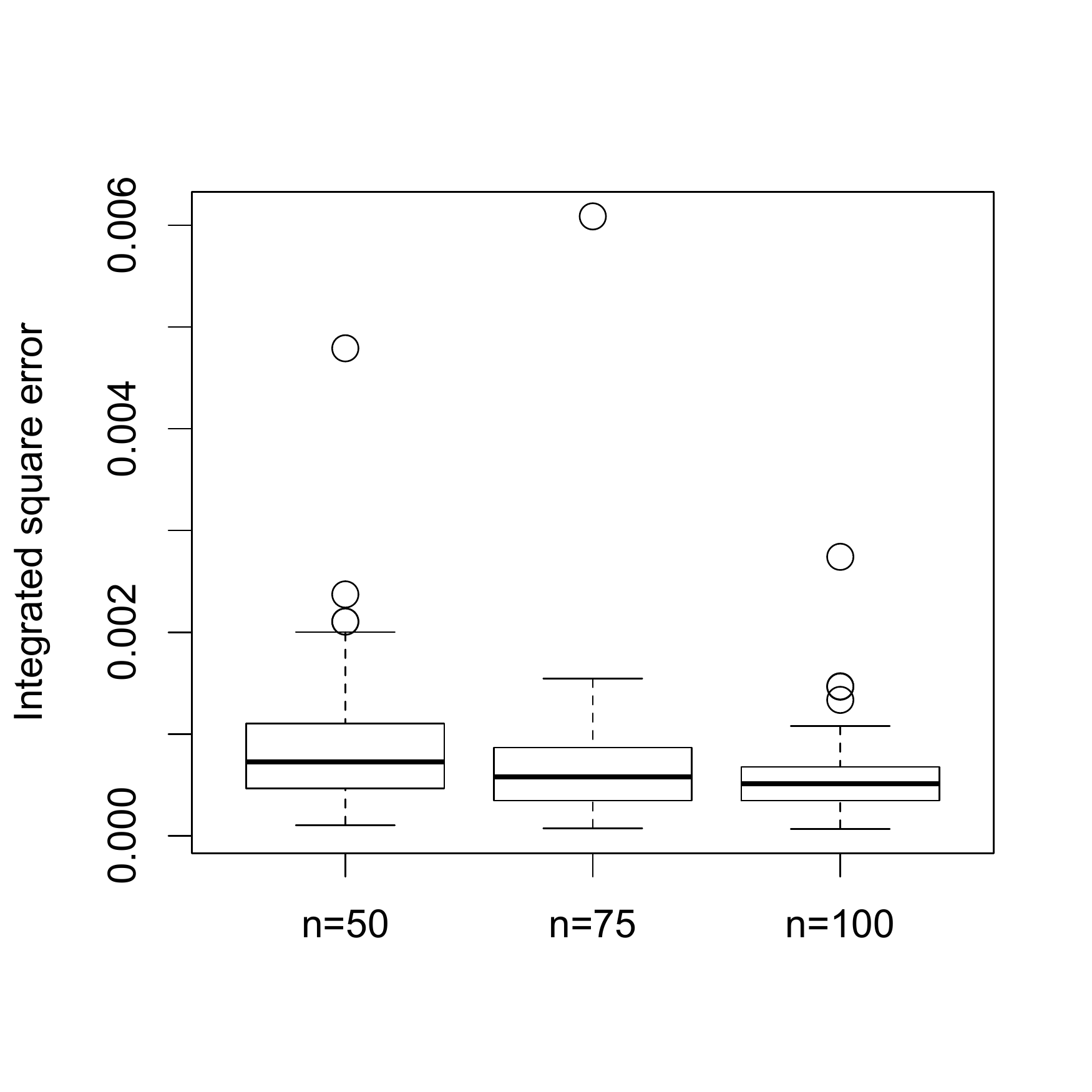}
	\qquad\qquad
	\includegraphics[height=5cm]{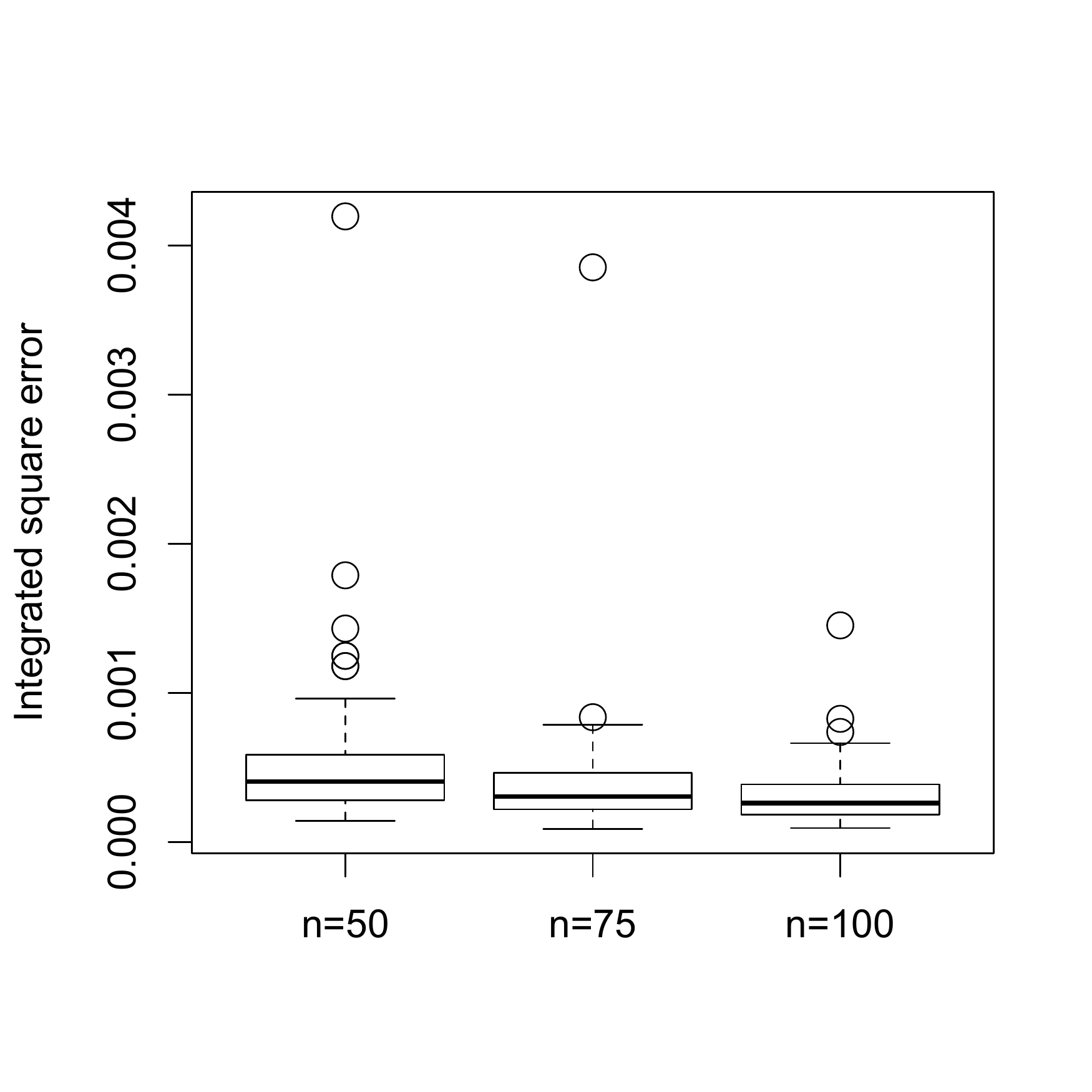}
	\includegraphics[height=5cm]{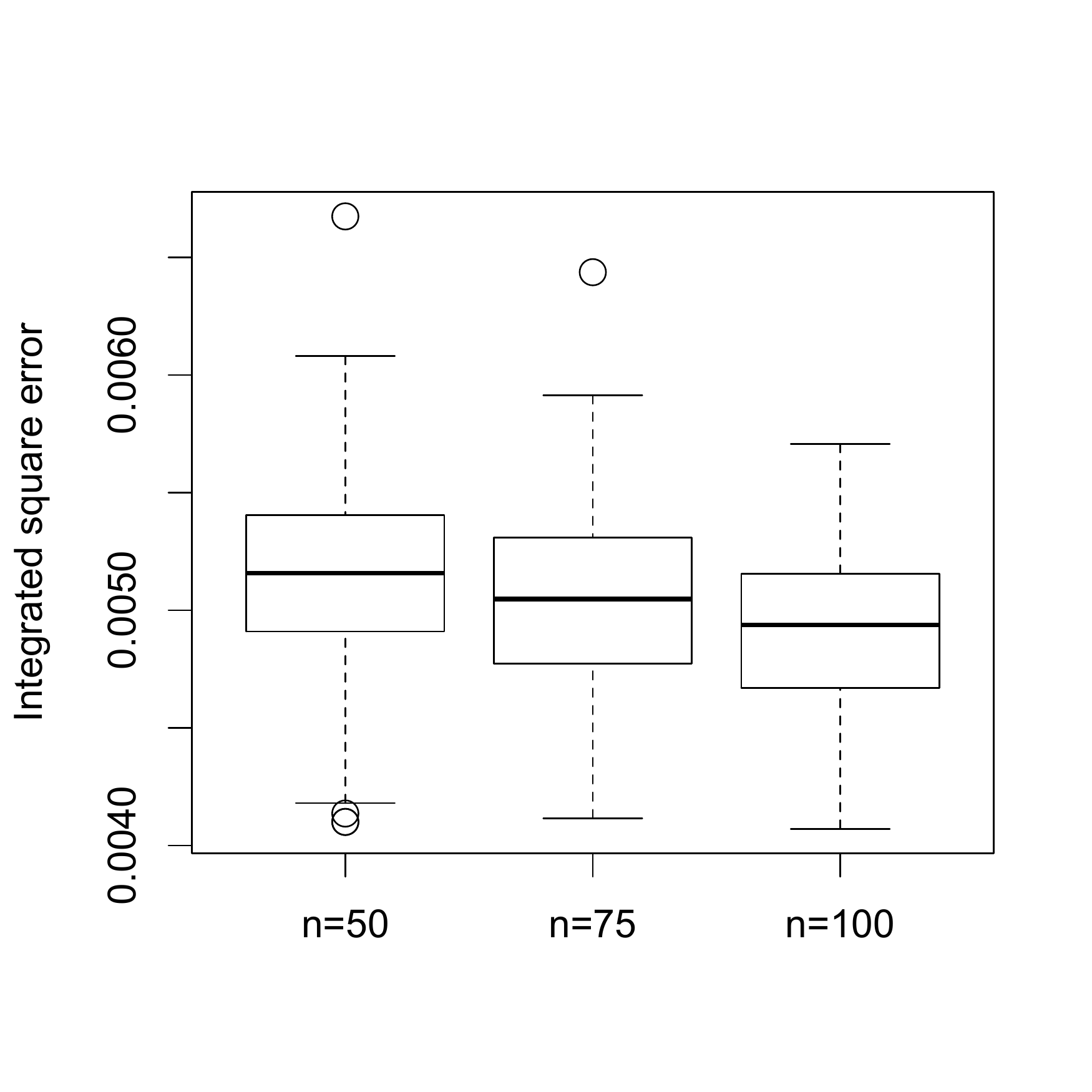}
	\qquad\qquad
	\includegraphics[height=5cm]{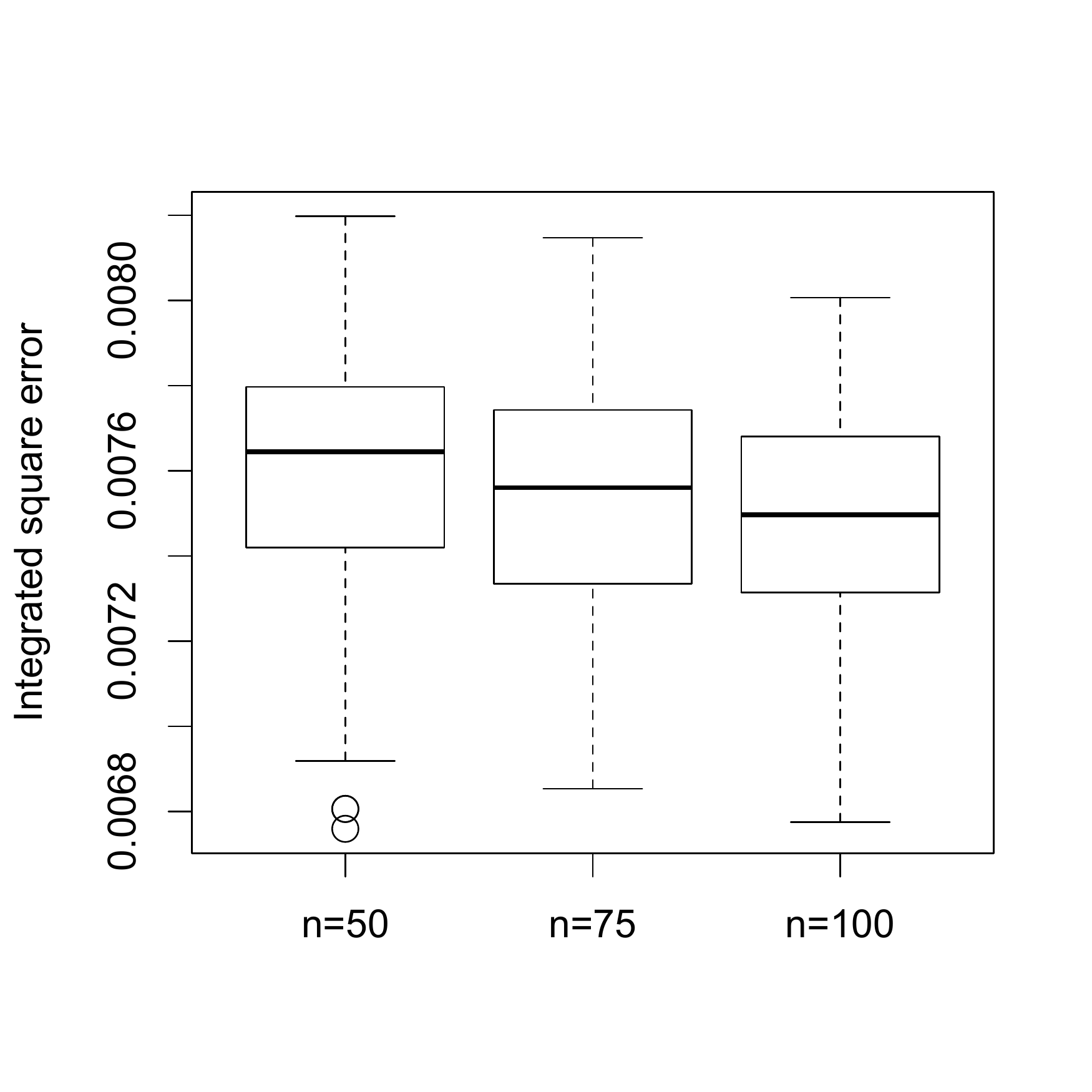}
	\caption{{The figure displays the} integrated square error on $100$ replicates between $t_m$ and its estimate $\widehat{t}_{m,n}$ from the observation of $n=50,\,75$ or $100$ random loss events and for $m=1$ (top left), $m=2$ (top right), $m=3$ (bottom left) and $m=4$ (bottom right).}
	\label{fig:t:ise}
\end{figure}

\begin{figure}[p]
	\centering
	\includegraphics[width=12cm,height=5.5cm]{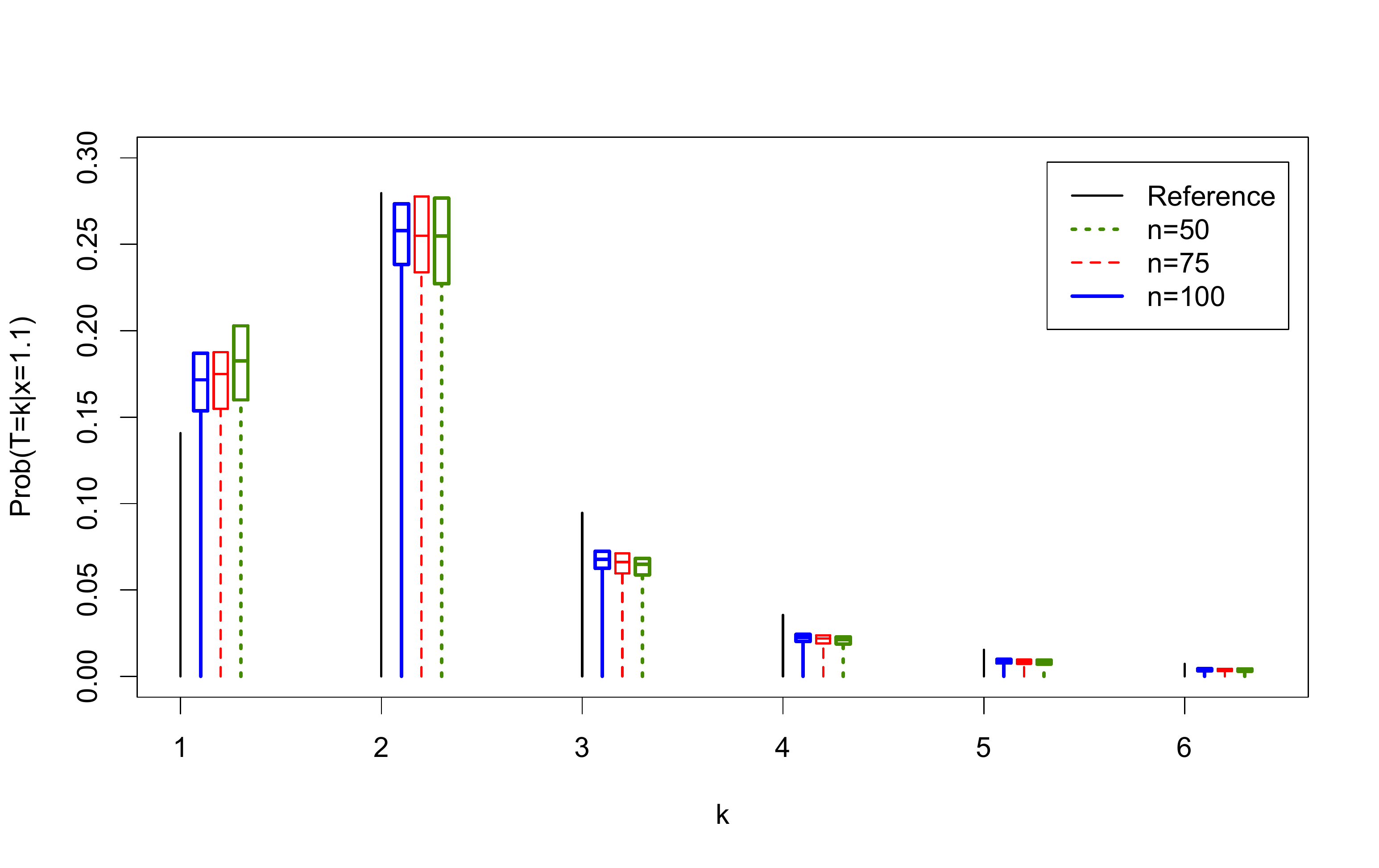}
	\caption{{The figure displays the} distribution of the hitting time of $\Gamma$ $t_m(x)$ for $x=1.1$ and $m=1,\,\dots,\,6$ and its estimates $\widehat{t}_{m,n}$ from the observation of $n=50,\,75$ or $100$ random loss events.}
	\label{fig:t}
\end{figure}

%\noindent
%\textbf{Romain Aza\"\i{}s}\\
%Inria Sophia Antipolis M\'editerran\'ee, Team Virtual Plants\\
%\url{romain.azais@gmail.com}
%
%\noindent
%\textbf{Alexandre Genadot}\\
%Laboratoire de Probabilit\'es et Mod\`eles Al\'eatoires, Universit\'e Pierre et Marie Curie, Paris 6\\
%\url{algenadot@gmail.com}
\end{document}